\documentclass[12pt]{article}
\include{epsf}
\usepackage{amsfonts}
\usepackage{graphics}
\usepackage{psfrag}
\usepackage{float}
\newtheorem{THEOREM}{Theorem}[section]
\newtheorem{PROPOSITION}{Proposition}[section]
\newtheorem{LEMMA}{Lemma}[section]

\newcommand{\qed}{\ \rule[-1pt]{4pt}{8pt}

                                         \vspace{2ex} }
\newenvironment{PROOF}{

                       \noindent{\bf Proof}.}{\qed}

\newcounter{labelflag} 

\newcommand{\Label}[1]{
                       \ifnum\thelabelflag=1
                          \ifmmode
                             \makebox[0in][l]{\qquad\fbox{\rm#1}}
                          \else
                             \marginpar{\vspace{0.7\baselineskip}
                                        \hspace{-1.1\textwidth}
                                        \fbox{\rm#1}}
                          \fi
                       \fi
                       \label{#1}
                      }
\newcommand{\be}{\begin{equation}}
\newcommand{\ee}{\end{equation}}
\newcommand{\eps}{\varepsilon}
\newcommand{\hm}{h_m}
\newcommand{\hmm}{h_{m+1}}
\newcommand{\ho}{h_{[0]}}
\newcommand{\hone}{h_{[1]}}
\newcommand{\hi}{h_{[i]}}
\newcommand{\hell}{h_{[\ell]}}
\newcommand{\hmone}{h_{[m+1]}}
\newcommand{\Me}{{\mathcal L}}
\newcommand{\Mo}{{\mathcal L}_{[0]}}

\newcommand{\Lo}{{\mathcal L}_0}
\newcommand{\Lm}{{\mathcal L}_{m}}
\newcommand{\Lmm}{{\mathcal L}_{m+1}}

\newcommand{\Ns}{{\mathrm N_{\mathrm s}}}
\newcommand{\Nf}{{\mathrm N_{\mathrm f}}}
\newcommand{\Nfs}{{\mathrm N}}
\newcommand{\R}{\mathbf{R}}
\newcommand{\T}{\mathrm T}
\newcommand{\N}{\mathrm N}
\newcommand{\hide}[1]{}
\newcommand{\abs}[1]{\left\vert#1\right\vert}
\newcommand{\Hmax}{H^{\rm max}}
\begin{document}

\begin{center}
{\large\textbf{
Analysis of the accuracy and convergence
of equation-free projection to a slow manifold
}}\\[1ex]
May 14, 2007

\vspace{2ex}
A.~Zagaris$^{1,2}$,
C.~W.~Gear$^{3,4}$,
T.~J.~Kaper$^5$,
I.~G.~Kevrekidis$^{3,6}$

$1$~{\small Department of Mathematics,
University of Amsterdam,
Amsterdam, The Netherlands.} \\
$2$~{\small Modeling, Analysis and Simulation,
Centrum voor Wiskunde en Informatica,
Amsterdam, The Netherlands.} \\
$3$~{\small Department of Chemical Engineering,
Princeton University,
Princeton, NJ 08544;} \\
$4$~{\small NEC Laboratories USA,
retired;} \\
$5$~{\small Department of Mathematics and Center for BioDynamics,
Boston University,
Boston, MA 02215;} \\
$6$~{\small Program in Applied and Computational Mathematics,
Princeton University,
Princeton, NJ 08544;}
\end{center}


\begin{abstract}
\noindent
In [C.~W.~Gear, T.~J.~Kaper, I.~G.~Kevrekidis, and A.~Zagaris,
Projecting to a Slow Manifold:
Singularly Perturbed Systems and Legacy Codes,
{\it SIAM J. Appl. Dyn. Syst.}
\textbf{4} (2005) 711--732],
we developed
a class of iterative algorithms
within the context
of equation-free methods
to approximate
low-dimensional,
attracting,
slow manifolds
in systems
of differential equations
with multiple time scales.
For user-specified values
of a finite number
of the observables,
the $m-$th member
of the class
of algorithms
($m = 0, 1, \ldots$)
finds iteratively
an approximation
of the appropriate zero
of the $(m+1)-$st time derivative
of the remaining variables
and
uses this root
to approximate the location
of the point
on the slow manifold
corresponding to these values
of the observables.
This article
is the first
of two articles
in which
the accuracy and convergence
of the iterative algorithms
are analyzed.
Here,
we work directly
with explicit fast--slow systems,
in which there is
an explicit small parameter,
$\eps$,
measuring the separation
of time scales.
We show that,
for each
$m = 0, 1, \ldots$,
the fixed point
of the iterative algorithm
approximates the slow manifold
up to and including
terms of ${\mathcal O}(\eps^m)$.
Moreover,
for each $m$,
we identify explicitly
the conditions
under which
the $m-$th iterative algorithm
converges to this fixed point.
Finally,
we show that
when
the iteration
is unstable
(or
converges slowly)
it may be stabilized
(or
its convergence
may be accelerated)
by application
of the Recursive Projection Method.
Alternatively,
the Newton--Krylov
Generalized Minimal Residual Method
may be used.
In the subsequent article,
we will consider
the accuracy and convergence
of the iterative algorithms
for a broader class
of systems---in which
there need not be
an explicit
small parameter---to which
the algorithms also apply.
\end{abstract}

\section{Introduction \label{s-intro}}
\setcounter{equation}{0}
The long-term dynamics
of many complex
chemical,
physical,
and biological systems
simplify
when a low-dimensional,
attracting,
invariant slow manifold
is present.
Such a slow manifold
attracts all nearby initial data
exponentially,
and
the reduced dynamics
on it
govern
the long term evolution
of the full system.
More specifically,
a slow manifold
is parametrized
by observables
which are typically
slow variables
or
functions of variables.
All nearby system trajectories
decompose naturally into
a fast component
that contracts exponentially
toward the slow manifold
and
a slow component
which obeys
the reduced system dynamics
on the manifold.
In this sense,
the fast variables
become slaved
to the observables,
and
knowledge of
the slow manifold
and
of the reduced dynamics
on it
suffices to determine
the full long-term system dynamics.

The identification and approximation
of slow manifolds
is usually achieved
by employing
a reduction method.
We briefly list
a number of these:
Intrinsic Low Dimensional Manifold
(ILDM),
Computational Singular Perturbation
(CSP),
Method of Invariant Manifold
(MIM),
Approximate Inertial Manifold approaches,
and
Fraser-Roussel iteration,
and we refer
the reader
to \cite{GKKZ-2005,KK-2002}
for a more extensive listing.

\subsection{A class of iterative algorithms  
based on the zero-derivative principle}
\label{subsec-1}

In~\cite{GKKZ-2005},
we developed a class
of iterative algorithms
to locate slow manifolds
for systems of
Ordinary Differential Equations
(ODEs)
of the form
\begin{eqnarray}
\begin{array}{rcccc}
    u' &=& p(u,v) ,
&\qquad&
    u \in \R^\Ns , \\
    v' &=& q(u,v) ,
&\qquad&
    v \in \R^\Nf ,
\end{array}
\Label{w-system}
\end{eqnarray}
where
$\Ns + \Nf \equiv \Nfs$.
We treated
the variables $u$
as the observables
(that is,
as parametrizing
the slow manifold
we are interested in),
and
we assumed that
there exists
an $\Ns-$dimensional,
attracting,
invariant,
slow manifold $\Me$,
which is given locally
by the graph
of a function
$v = v(u)$.
However,
we emphasize
that we did not need
explicit knowledge
of which variables
are fast
and
which are slow,
only that
the variables $u$
suffice to parametrize $\Me$.

To leading order,
the location 
of a slow manifold $\Me$
is obtained by setting $v'=0$,
{\it i.e.,} by solving
$q(u,v)=0$ for $v$.
Of course,
the manifold defined
by this equation
is in general not
an invariant slow manifold
under the flow of the full system
(\ref{w-system}).
This is only approximately true,
since higher-order derivatives
with respect to the (fast) time $t$
are, in general, large on it.
If one requires that
$v''$ vanishes,
then the solutions
with initial conditions at the points
defined by this condition
depend only on the slow time
to one order higher,
as $v'$ also remains bounded
in the vicinity of this manifold.
Similarly,
demanding that successively higher-order
time derivatives vanish,
we obtain manifolds
where all time derivatives
of lower order remain bounded.
The solutions with these initial conditions
depend only on the slow time
to successively higher order
and thus approximate,
also to successively higher order,
solutions on the slow manifold.
In other words,
demanding that time derivatives
of successively higher order vanish,
we filter out the fast dynamics
of the solutions
to successively higher orders.
In this manner, the approximation
of the slow manifold $\Me$ 
is improved successively,
as well.
This idea 
may be traced back
at least to the work 
of Kreiss~\cite{BK-1982, K-1979, K-1985},
who studied systems 
with rapid oscillations
(asymptotically large frequencies)
and introduced 
the bounded derivative principle
to find approximations of slow manifolds
as the sets of points
at which the derivatives 
are bounded (not large).
The requirement here
that the derivatives 
with respect 
to the (fast) time $t$ vanish
is the analog 
for systems (\ref{w-system})
with asymptotically stable 
slow manifolds.
A similar idea was introduced independently
by Lorenz in~\cite{L-1980},
where he used
a simple functional iteration scheme
to approximate the zero
of the first derivative,
then used the converged value
of this scheme
to initialize a similar scheme
that approximates the zero
of the second derivative,
and so on until successive zeroes
were found to be virtually identical.
See also \cite{CHL-1995} and \cite{G-1999}
for other works in which
a similar condition is employed.

The elements 
of the class
of iterative algorithms
introduced in \cite{GKKZ-2005}
are indexed by $m=0,1,\ldots$.
The $m-$th algorithm
is designed to locate,
for any fixed value
of the observable $u_0$,
an appropriate solution,
$v=v_m(u_0)$,
of the $(m+1)-$st
derivative condition
\be
  \left( \frac{d^{m+1} v}{dt^{m+1}} \right)(u_0,v)
=
  0 .
\Label{deriv-cond-w}
\ee
Here,
the time derivatives
are evaluated along solutions
of (\ref{w-system}).
In general,
since condition (\ref{deriv-cond-w})
constitutes a system
of $\Nf$ nonlinear algebraic equations,
the solution $v_m(u_0)$
cannot be computed explicitly.
Also,
the explicit form of (\ref{w-system}),
and thus also an analytic formula
for the $(m+1)-$st time derivative
in Eq.~(\ref{deriv-cond-w}),
may be unavailable
(\emph{e.g.}, in Equation-Free 
or legacy code applications).
In this case,
a numerical approximation for it
has to be used.
The $m$-th algorithm in the class
generates an approximation
$v^\#_m$ of $v_m(u_0)$,
rather than
$v_m(u_0)$ itself,
using either an analytic formula
for the time derivative
or a finite difference approximation for it.
In either case,
the approximation $v^\#_m$ to $v_m(u_0)$
is determined through an explicit
functional iteration scheme,
which we now introduce.

The $m=0$ algorithm
is defined by the map
${\tilde F}_0 : \R^\Nf \to \R^\Nf$
\[
  {\tilde F}_0(v)
=
  v
+
  H \left( \frac{dv}{dt} \right)(u_0 , v) ,
\]
where
$H$, which we 
label as the iterative step size, is
an arbitrary positive number
whose magnitude
we fix below
for stability reasons.
We initialize
the iteration
with some value
$v^{(1)}$
and
generate the sequence
\[
  \left\{\left.
  v^{(r+1)} \equiv {\tilde F}_0( v^{(r)} )
  \right\vert
  r = 1,2, \ldots
  \right\} .
\]
The functional iteration
is terminated
when
$\| v^{(r+1)} - v^{(r)} \|
<
\mathrm{TOL}_0$,
for some $r \ge 1$
and
a prescribed
tolerance $\mathrm{TOL}_0$.
The output
of this zeroth algorithm
is
the last member,
$v^\#_0$,
of the sequence
$\{ v^{(r+1)} \}$.

Next,
the $m=1$ algorithm
is defined by the map
${\tilde F}_1 : \R^\Nf \to \R^\Nf$,
\[
  {\tilde F}_1(v)
=
  v
-
  H^2 \left( \frac{d^2v}{dt^2} \right)(u_0 , v) ,
\]
initialized with some value $v^{(1)}$.
It generates the sequence
\[
  \left\{\left.
  v^{(r+1)} \equiv {\tilde F}_1( v^{(r)} )
  \right\vert
  r = 1, 2, \ldots
  \right\} 
\]
and the functional iteration
is terminated
when
$\| v^{(r+1)} - v^{(r)} \|
<
\mathrm{TOL}_1$,
for some $r \ge 1$
and
for a prescribed
tolerance $\mathrm{TOL}_1$.
The output
of this first algorithm
is the last member,
$v^\#_1$,
of the sequence
$\{ v^{(r+1)} \}$.

The algorithm with general $m$
is defined by the map
${\tilde F}_m : \R^\Nf \to \R^\Nf$,
\be
  {\tilde F}_m(v)
=
  v
-
  (-H)^{m+1} \left( \frac{d^{m+1}v}{dt^{m+1}} \right)(u_0 , v) ,
\Label{tilde-F(m)-def-raw}
\ee
seeded with some value $v^{(1)}$.
It generates the sequence
\[
    \left\{\left.
    v^{(r+1)} \equiv {\tilde F}_m( v^{(r)} )
    \right\vert
    r = 1, 2, \ldots
    \right\} .
\]
Here also,
one prescribes
a tolerance $\mathrm{TOL}_m$
and
terminates the iteration procedure
when
$\| v^{(r+1)} - v^{(r)} \|
<
\mathrm{TOL}_m$
for some $r \ge 1$.
The output
of this $m-$th algorithm
is the last member
of the sequence
$\{ v^{(r+1)} \}$,
denoted by
$v^\#_m$.

As we show in this article,
not only is the point
$(u_0 , v^\#_m)$ of interest
for each individual $m$
because it \emph{approximates}
$(u_0 , v(u_0))$,
but the entire sequence
$\{ (u_0 , v^\#_m) \}_m$
is also of interest
because
it \emph{converges}
to $(u_0 , v(u_0))$
with a suitably convergent
sequence $\{ \mathrm{TOL}_m \}$.
Hence, the latter point
can be approximated
arbitrarily well
by members of that sequence,
and the class of algorithms
may be used 
as an integrated sequence
of algorithms
in which
the output $v^\#_m$
of the $m-$th algorithm
can be used to initialize
the $(m+1)-$st algorithm.
Of course,
other initializations
are also possible,
and we have carried
out the analysis here
in a manner
that is independent
of which choice one makes.

This class of iterative algorithms
was applied in \cite{GKKZ-2005}
to three examples:
the two-dimensional
Michaelis--Menten mechanism
for which
the one-dimensional slow manifold
can be computed analytically
to arbitrary precision,
a five-dimensional nonlinear system
with an explicitly computable
two-dimensional slow manifold,
and a seven-dimensional
hydrogen-oxygen system
with quadratic nonlinearities
for which the manifold
is not known explicitly.
In the context 
of these three examples,
we found that,
for all of the values of $m$ 
that we worked with,
the $m$-th algorithm 
converged at an exponential rate
to a fixed point.
Moreover,
in the two examples
where the slow manifold
can be computed,
we also found that,
for each algorithm,
this fixed point is
very close
to the actual point
on the slow manifold.
In addition
to showing
the $m$-th algorithm
converged for each $m$ 
that we worked with,
we also showed 
that the class of algorithms
may be used 
in the integrated manner
stated above.
The closeness 
of the approximation
to $(u_0,v(u_0))$
improved as we increased 
the order $m$ of the algorithm used.

More recently,
van Leemput \emph{et al.} \cite{vLVR-2005}
employed the first ($m=0$) algorithm
in the class
to initialize Lattice Boltzmann Models (LBM)
from sets of macroscopic data
in a way that eliminates
the stiff dynamics
triggered by a bad initialization.
They showed that
the algorithm they derived
converges unconditionally
to a fixed point
close to a slow manifold,
and
they used the algorithm
to couple a LBM
to a reaction-diffusion equation
along the interface
with good results
\cite{vLVR-2006}.

Our motivation
for introducing
this class
of iterative algorithms
in \cite{GKKZ-2005}
was two-fold.
First,
we wanted a method
that can be implemented
in the context
of legacy codes.
In other words,
we wanted
this reduction method
to be implementable
even when one
has no explicit form
for the components
$p$ and $q$
of the vector field,
but only
a black-box integrator
(timestepper).
This feature
renders the method
``equation-free'' \cite{KGHKOT-2003}
and
makes its implementation
possible in these settings.
Second,
it was essential
for us
that they preserve
the user-specified value
of the observables,
say $u = u_0$,
at each iteration.
In this way,
the output
of the algorithm
is an approximation
of the point
$(u_0 , v(u_0))$
on the manifold $\Me$
of that same value
$u_0$
of the observables.
Also,
in this way,
the `lifting' step
in projective integration
of \cite{GK-2004}
is naturally facilitated.

It is worth noting that
one really 
only needs to require
that the time derivatives 
are sufficiently small, 
although we work 
with the zero-derivative condition
(\ref{deriv-cond-w})
for definiteness.

\subsection{Iterative algorithms
based on the zero-derivative principle
for explicit fast--slow systems}
\label{subsec-2}

A central assumption
that we made
in \cite{GKKZ-2005}
is that
we work with
systems~(\ref{w-system})
for which
there exists
a smooth and invertible
coordinate change
\be
    z = z(w)
\quad\mbox{with inverse}\quad
    w = w(z) ,
\Label{coc}
\ee
where
$w = (u,v)$
and
$z = (x,y)$,
which puts the system~(\ref{w-system})
into the explicit fast--slow form
\begin{eqnarray}
\begin{array}{rcccc}
    x' &=& f(x,y,\eps) ,
&\qquad&
    x \in \R^\Ns , \\
    \eps y' &=& g(x,y,\eps) ,
&\qquad&
    y \in \R^\Nf .
\end{array}
\Label{z-system}
\end{eqnarray}
We emphasize that,
in general,
we have
no knowledge whatsoever
of the transformation
that puts system~(\ref{w-system})
into an explicit fast--slow form.
Here,
$f$ and $g$
are smooth functions
of their arguments,
the manifold $\Me$
is transformed smoothly,
and
$\mathrm{det}(D_y g)_0(z)
\equiv
\mathrm{det}(D_y g(z,0))
\ne
0$
on the manifold
$\Mo = \{ z \vert g(z,0) = 0 \}$
(on which
the dynamics reduce
for $\eps = 0$),
see also~\cite{GKKZ-2005}.

Due to the above assumption,
it turns out to be natural
to split the analysis
of the accuracy and convergence
of the functional iteration
into two parts.
In the first part,
which we present
in this article,
we work directly
on systems that are already
in explicit fast--slow form (\ref{z-system}).
In the context of these systems,
the accuracy and convergence analysis
may be carried out completely
in terms of the small
parameter $\eps$.
The system geometry --
the slow manifold
and the fast fibers transverse to ${\cal L}$ --
makes the convergence analysis
especially transparent.
Then,
in the second part,
we work with the more general systems
(\ref{w-system}).
For these,
the accuracy analysis proceeds
along similar lines
as that for this first part,
with the same type of result
as Theorem~\ref{t-main1} below.
However,
the convergence analysis
is considerably more involved
than that for explicit fast--slow systems.
For these general systems,
one must analyze
a series of different scenarios
depending on the relative orientations
of (i) the tangent space to ${\cal L}$,
(ii) the tangent spaces to the fast fibers
at their base points on ${\cal L}$, and
(iii) the hyperplane of the observables $u$.
Moreover, all of the analysis
must be carried out
through the lens
of the coordinate change (\ref{coc})
and its inverse,
so that it
is less transparent
than it is in part one.
Part two will be presented 
as a subsequent article.

As applied specifically
to explicit fast--slow systems~(\ref{z-system}),
the $m-$th iterative algorithm
(\ref{tilde-F(m)-def-raw})
is based on
the $(m+1)-$st
derivative condition,
\be
  \left( \frac{d^{m+1} y}{dt^{m+1}} \right)(x_0 , y)
=
    0 .
\Label{deriv-cond}
\ee
In particular,
for each $m$
and
for any arbitrary,
but fixed,
value of the observable
$x_0 \in K$,
one makes
an initial guess
for $h(x_0)$ and
uses the $m$-th iterative algorithm
to approximate
the appropriate zero
of this $(m+1)-$st derivative,
where the end (\emph{converged}) result
of the iteration
is the improved approximation
of $h(x_0)$.

For each
$m = 0, 1, \ldots$,
the $m-$th iterative algorithm
is defined
by the map
$F_m : \R^\Nf \to \R^\Nf$,
\be
  F_m(y)
=
  y
-
  (-H)^{m+1} \left( \frac{d^{m+1}y}{dt^{m+1}} \right)(x_0 , y) ,
\Label{F(m)-def-raw}
\ee
where $H$ is
an arbitrary positive number
whose magnitude is
${\mathcal O}(\eps)$
for stability reasons.
We seed with some value $y^{(1)}$
and
generate the sequence
\be
    \left\{\left.
    y^{(r+1)} \equiv F_m( y^{(r)} )
    \right\vert
    r = 1, 2, \ldots
    \right\} .
\Label{seq(m)}
\ee
Here also,
one prescribes
a tolerance $\mathrm{TOL}_m$
and
terminates the iteration procedure
when
$\| y^{(r+1)} - y^{(r)} \|
<
\mathrm{TOL}_m$
for some $r \ge 1$.
The output
of this $m-$th algorithm
is the last member
of the sequence
$\{ y^{(r+1)} \}$,
denoted by
$y^\#_m$.

\subsection{Statement of the main results}
\label{subsec-3}

In this article,
we first examine
the $m$-th iterative algorithm
in which an analytical formula
for the $(m+1)-$st derivative
is used,
and we prove that it 
has a fixed point
$y=h_m(x_0)$,
which is ${\cal O}(\eps^{m+1})$ close
to the corresponding point $h(x_0)$
on the invariant manifold ${\cal L}$,
for each $m=0,1,\ldots$.
See Theorem~\ref{t-main1} below.

Second, we determine
the conditions on $(D_y g)_0$
under which 
the $m$-th iterative algorithm converges
to this fixed point,
again with an analytical formula
for the $(m+1)-$st derivative.
In particular, for $m=0$,
the iteration converges
for all systems (\ref{z-system})
for which $(D_y g)_0$
is uniformly Hurwitz on $\Mo$
and
provided that
the iterative step size $H$ is small enough.
For each $m \ge 1$,
convergence of the algorithm
imposes more stringent conditions
on $H$ and
on the spectrum
of $(D_y g)_0$.
In particular,
if
$\sigma((D_y g)_0)$ is contained
in certain sets in the complex plane,
which we identify completely,
then
the iteration converges
for small enough values
of the iterative step size $H$,
see Theorem~\ref{t-main2}.
These sets
do not cover
the entire half-plane,
and thus
complex eigenvalues can,
in general,
make the algorithm divergent.

Third, we show explicitly
how the Recursive Projection Method (RPM)
of Shroff and Keller \cite{SK-1993}
stabilizes the functional iteration
for each $m \ge 1$
in those regimes where the iteration is unstable.
This stabilization result is useful
for practical implementation
in the equation-free context; and,
the RPM may also be used to accelerate convergence
in those regimes
in which the iterations converge slowly.
Alternatively,
the Newton--Krylov
Generalized Minimal Residual Method
(NK-GMRES \cite{K-1995})
may be used to achieve
this stabilization.

Fourth, we analyze
the influence of the tolerance,
or stopping criterion,
used to terminate the functional iteration.
We show that,
when the tolerance
$\mathrm{TOL}_m$
for the $m-$th algorithm
is set to
${\mathcal O}(\eps^{m+1})$,
the output $y^\#_m$
also satisfies
the asymptotic estimate
$\| y^\#_m - h(x_0) \|
=
{\mathcal O}(\eps^{m+1})$.

Finally, we extend
the accuracy and convergence analyses
to the case where
a forward difference approximation
of the $(m+1)-$st derivative
is used
in the iteration,
instead of the analytical formula.
As to the accuracy,
we find that
the $m$-th iterative algorithm
also has a fixed point
$y=\hat{h}_m(x_0)$
which is ${\cal O}(\eps^{m+1})$ close
to $h(x_0)$,
so that
the iteration in this case
is as accurate asymptotically
as the iteration 
with the analytical formula.
Then, as to the stability,
we find that
the $m$-th iterative algorithm 
with a forward difference approximation
of the $(m+1)-$st derivative
converges unconditionally
for $m=0$.
Moreover,
for $m=1,2,\ldots$,
the convergence 
is for a continuum of values
of the iterative step size $H$
and
without further restrictions
on $(D_y g)_0$,
other than that
it is uniformly Hurwitz
on $\Mo$,
see Theorem~\ref{t-main3}.
These advantages stem from the use
of a forward difference approximation,
and we will show in a future work
that the use of implicitly defined
maps $F_m$ yields similar advantages.

Throughout this article,
we shall refer to
some basic facts about
the $\Ns-$dimensional,
slow, invariant,
and
normally attracting
manifold $\Me$.
As stated above,
$\Me$
is the graph
of a function $h$,
\be
  \Me
=
  \left\{
  (x,y)
  \left\vert
  x \in K,\,
  y = h(x)
  \right.
  \right\} ,
\Label{Me}
\ee
for some set $K$.
Here, the function
$h : K \to \R^\Nf$
satisfies
the invariance equation
\be
  g(x,h(x),\eps)
-
  \eps Dh(x) f(x , h(x),\eps)
=
  0 ,
\Label{inv-eq}
\ee
and it is
${\cal O}(\eps)$ close
to the critical manifold,
which is the graph of $h_0(x)$,
uniformly for $x \in K$.

It is insightful
to recast
this invariance equation
in the form
\be
  \left( - Dh(x) , I_\Nf \right)
  G(x , h(x),\eps)
=
  0 ,
\quad \mbox{where} \quad
  G
\equiv
  \left(
  \begin{array}{c}
  \eps f \\ g
  \end{array}
  \right) ,
\Label{inv-eq-alt}
\ee
which reveals
a clear geometric interpretation.
Since $\Me$ corresponds
to the zero level set
of the function
$- h(x) + y$
by Eq.~(\ref{Me}),
the rows of
the $\Nf \times \Nfs$
gradient matrix
$(- Dh(x) , I_\Nf )$
form a basis
for $\N_z \Me$,
the space normal
to the slow manifold
at the point
$z = (x , h(x)) \in \Me$.
Thus,
Eq.~(\ref{inv-eq-alt}) states
that the vector field $G$
is perpendicular to this space
and
hence contained
in the space
tangent to
the slow manifold,
$\T_z \Me$.

\section{Existence of a fixed point $\hm(x_0)$
and its proximity to $h(x_0)$ \label{s-stage1}}
\setcounter{equation}{0}
We rewrite
the map $F_m$,
given in Eq.~(\ref{F(m)-def-raw}),
as
\be
  F_m(y)
=
  y
-
  L_m(x_0 , y) ,
\Label{F(m)-def}
\ee
where
the function
$L_m : \R^\Nfs \to \R^\Nf$
is given by
\be
  L_m(z)
\equiv
  (-H)^{m+1}
  \left( \frac{d^{m+1} y}{dt^{m+1}} \right)(z) ,
\quad\mbox{for any}\:
  m = 0,1,\ldots\, ,
\Label{L(m)-def}
\ee
where
$z = (x_0 , y)$.
The fixed points,
$y = \hm(x_0)$,
of $F_m$
are determined
by the equation
\[
    L_m(x_0 , \hm(x_0))
=
    0 ,
\]
that is,
by the $(m+1)-$st
derivative condition~(\ref{deriv-cond}).
The desired results
on the existence
of the
fixed point $\hm(x_0)$
and on its proximity
to $h(x_0)$
are then immediately at hand
from the following theorem:
\begin{THEOREM}
\label{t-main1}
For each $m=0,1,\ldots$,
the $(m+1)-$st
derivative condition~(\ref{deriv-cond}),
\be
  L_m(x,y)
\equiv
  (-H)^{m+1}
  \left( \frac{d^{m+1} y}{dt^{m+1}} \right)(x,y)
=
  0 ,
\Label{cond}
\ee
can be solved for $y$
to yield
an $\Ns-$dimensional manifold $\Lm$
which is the graph
of a function
$\hm : K \to \R^\Nf$
over $x$.
Moreover,
the asymptotic expansions
of $\hm$ and $h$
agree up to and including
terms of ${\mathcal O}(\eps^m)$,
\[
  \hm(\cdot)
=
  \sum_{i=0} \eps^i h_{m,i}(\cdot)
=
  \sum_{i=1}^m \eps^i \hi(\cdot)
+
  {\mathcal O}(\eps^{m+1}) .
\]
\end{THEOREM}
This theorem guarantees that,
for each $x_0 \in K$,
there exists an isolated fixed point
$y = \hm(x_0)$
of the functional iteration algorithm.
Moreover,
this fixed point
varies smoothly
with $x_0$,
and the approximation
$(x_0 , \hm(x_0))$
of the point
$(x_0 , h(x_0))$
on the actual
invariant slow manifold
is valid
up to ${\cal O}(\eps^{m+1})$.

The remainder of this section
is devoted to the proof
of this theorem.
We prove it
for $m=0$ and $m=1$
in Sections~\ref{ss-m=0} and \ref{ss-m=1},
respectively.
Then, in Section~\ref{ss-ind},
we use induction
to prove the theorem
for general $m$.

\subsection{Proof of Theorem~\ref{t-main1} for $m=0$ \label{ss-m=0}}
We show,
for each
$x \in K$,
that $L_0(z)$
has a root
$y
=
h_0(x)$,
that $h_0$
lies ${\mathcal O}(\eps)$
close to $\ho(x)$,
the corresponding point
on the critical manifold,
and that
the graph
of the function $h_0$
over $K$
forms a manifold.

For $m=0$,
definition~(\ref{L(m)-def}),
the chain rule,
and
the ODEs~(\ref{z-system})
yield
\begin{eqnarray}
  L_0
=
  -H y'
=
  -\eps^{-1} H g .
\Label{L(0)-raw}
\end{eqnarray}
Substituting the asymptotic expansion
$y
=
h_0(x)
=
\sum_{i=0} \eps^i h_{0,i}(x)$
into this formula
and
combining it
with the condition
$L_0 = 0$,
we find that,
to leading order,
\begin{eqnarray*}
  g(x, h_{0,0}(x), 0)
=
  0 ,
\quad
\end{eqnarray*}
where
we have removed
the ${\mathcal O}(1)$,
nonzero,
scalar quantity $-H/\eps$.
In comparison,
the invariance equation~(\ref{inv-eq})
yields
\begin{equation}
  g\left(x, \ho(x), 0\right)
=
  0 ,
\Label{inveq-0}
\end{equation}
to leading order,
see Eq.~(\ref{h(0)})
in Appendix A.
Thus
$h_{0,0}$ can be chosen
to be equal to $\ho$,
and
$L_0(z)$ has a root
that is ${\mathcal O}(\eps)-$close
to $y = h(x)$.

It remains
to show
that the graph
of the function $h_0$
is an $\Ns-$dimensional manifold
$\Lo$.
Using Eq.~(\ref{L(0)-raw}),
we calculate
\[
  \left(D_y L_0\right)
=
  - \eps^{-1} H \left(D_yg\right) ,
\]
where all quantities
are evaluated at
$(x, h_0(x), \eps)$.
Moreover,
\[
  \left(D_y L_0\right)(x, h_0(x))
=
  - \eps^{-1} H \left(D_yg\right)_0
+
  {\mathcal O}(\eps) ,
\]
with
$(\cdot)_0
=
(\cdot)(x, h_{0,0}(x), 0)
=
(\cdot)(x, \ho(x), 0)$,
since
$h_{0,0} = \ho$.
Thus,
the Jacobian
$(D_y L_0)(x, h_0(x))$
is non-singular
for $0 < \eps \ll 1$,
because
$H={\mathcal O}(\eps)$
by assumption
and because
$\mathrm{det}(D_y g)_0
\ne
0$,
see the Introduction.
Therefore,
we have
\[
  \mathrm{det}
  \left(D_y L_0\right)(x, h_0(x))
\ne
  0 ,
\quad \mbox{for all} \
  x \in K,
\]
and hence
$\Lo$
is a manifold
by the Implicit Function Theorem
and
\cite[Theorem~1.13]{O-1986}.
This completes
the proof
of the theorem
for the case $m=0$.

\subsection{The proof of Theorem~\ref{t-main1} for $m=1$ \label{ss-m=1}}
In this section,
we treat
the $m=1$ case.
Technically speaking,
one may proceed directly
from the $m=0$ case
to the induction step
for general $m$.
Nevertheless,
we find it useful
to present a concrete instance
and a preview of the general case,
and hence we give a brief analysis
of the $m=1$ case here.

We calculate
\[
  L_1
=
  (-H)^2 y''
=
  -H (-H y')'
=
  -H L_0'
=
  - \eps^{-1} H ( D_z L_0 ) G .
\]
Using the ODEs~(\ref{z-system})
and
Eq.~(\ref{L(0)-raw}),
we rewrite this as
\begin{eqnarray}
  L_1
=
  \left( -\eps^{-1} H \right)^2
  \left[
  \eps (D_xg) f
+
  (D_yg) g
  \right] .
\Label{L(1)-raw}
\end{eqnarray}
We recall
that the solution
is denoted by
$y
=
h_1(x)$
and that we write
its asymptotic expansion as
$h_1(x)
=
\sum_{i=0} \eps^i h_{1,i}(x)$.
Substituting this expansion
into Eq.~(\ref{L(1)-raw})
and
recalling that
$H = {\mathcal O}(\eps)$,
we obtain at ${\mathcal O}(1)$
\[
  L_1
=
  (-\eps^{-1} H)^2
  \left(D_yg\right)_0 g_0
+
  {\mathcal O}(\eps),
\]
where
$(\cdot)_0
=
(\cdot)(x, h_{1,0}(x), 0)$.
Hence,
$y = \ho(x)$
is a root of $L_1$
to leading order
by Eq.~(\ref{inveq-0})
and
$\mathrm{det}(D_y g)_0
\ne
0$,
and therefore
$h_{1,0}$ can be selected
to be equal to $\ho$.

At ${\mathcal O}(\eps)$,
we obtain
\begin{equation}
(-\eps^{-1} H)^2  (D_y g)_0
\left[
  (D_y g)_0^{-1} (D_xg)_0 f_0
+
  (D_y g)_0 h_{1,1}
+
  (D_\eps g)_0
\right]
=
  0 ,
\Label{L(1)=0}
\end{equation}
where we used
the expansion
\[
  g(x , h_1 , \eps)
=
  g_0
+
  \eps
\left[
  (D_y g)_0 h_{1,1}
+
  (D_\eps g)_0
\right]
+
  {\mathcal O}(\eps^2)
\]
and that
$g_0
=
g(x , h_{1,0}, 0)
=
g(x , \ho, 0)$.
Differentiating
both members of
the identity
$g(x , \ho(x) , 0) = 0$
with respect to $x$,
we obtain
\[
  (D_xg)_0
+
  (D_yg)_0 (D\ho)
=
  0 ,
\]
whence
$(D_yg)_0^{-1} (D_xg)_0
=
- D\ho$.
Removing
the invertible prefactor
$(-H/\eps)^2 (D_yg)_0$,
we find that
Eq.~(\ref{L(1)=0})
becomes
\[
  -(D\ho) f_0
+
  (D_y g)_0 h_{1,1}
+
  (D_\eps g)_0
=
  0 .
\]
This equation
is identical to
Eq.~(\ref{h(1)})
in Appendix A,
and thus
$h_{1,1}=\hone$.
Hence, we have shown
that the asymptotic expansion
of $h_1(x)$
agrees with that of $h(x)$
up to and including
terms of ${\mathcal O} (\eps)$,
as claimed for $m=1$.

Finally,
the graph
of the function $h_1$
forms an
$\Ns-$dimensional manifold
${\mathcal L}_1$.
This may be shown
in a manner
similar to that
used above
for $\Lo$
in the case
$m=0$.
This completes the proof for $m=1$.

\subsection{The induction step: the proof of Theorem~\ref{t-main1} for general $m$ \label{ss-ind}}
In this section,
we prove
the induction step
that establishes
Theorem~\ref{t-main1}
for all $m$.
We assume that
the conclusion
of Theorem~\ref{t-main1}
is true for $m$
and show
that it also holds
for $m+1$,
{\it i.e.,}
that the condition
\begin{equation}
    \left[(D_z L_m)(x, y)\right]
    G(x, y, \eps)
=
    0
\Label{cond-alt-m+1}
\end{equation}
can be solved for $y$
to yield
$y
=
\hmm(x)$,
where
\[
    \hmm(\cdot)
=
    \sum_{i=0}^{m+1} \eps^i \hi(\cdot) + {\mathcal O}(\eps^{m+2}) .
\]
To begin with,
we recast
the $(m+1)-$st
derivative condition
Eq.~(\ref{cond})
in a form that
is reminiscent
of the invariance equation,
Eq.~(\ref{inv-eq-alt}).
Let $m \ge 0$
be arbitrary but fixed.
It follows
from definition~(\ref{L(m)-def}),
Eq.~(\ref{inv-eq-alt}),
and
Eq.~(\ref{z-system})
that
\be
    L_m
=
    -H \frac{d}{dt}\left( (-H)^m \frac{d^m y}{dt^m} \right)
=
    -H \frac{dL_{m-1}}{dt}
=
    - \eps^{-1} H ( D_z L_{m-1} ) G .
\Label{L(m+1)=DL(m)G}
\ee
Therefore,
the $(m+1)-$st derivative
condition~(\ref{cond})
can be rewritten
in the desired form
as
\begin{equation}
    \left( D_z L_{m-1} \right) G = 0 ,
\Label{cond-alt}
\end{equation}
where
we have removed
the ${\mathcal O}(1)$,
nonzero, scalar quantity $-H/\eps$.

The induction step
will be now be established
using a bootstrapping approach.
First,
we consider
a modified version
of Eq.~(\ref{cond-alt-m+1}),
namely the condition
\begin{equation}
    \left[(D_z L_m)(x, \hm(x))\right] G(x, y, \eps)
=
    0 ,
\Label{mod-cond}
\end{equation}
in which
the matrix
$D_z L_m$
is evaluated on $\Lm$
(already determined
at the $m-$th iteration)
instead of
on the as-yet
unknown $\Lmm$.
This equation
is easier
to solve
for the unknown $y$,
since
$y$ appears only in $G$.
We now show that
the solution
$y = \tilde{h}_{m+1}(x)$
of this condition
approximates $h$
up to and including
${\mathcal O}(\eps^{m+1})$ terms.
\begin{LEMMA}
\label{l-KGKvsCSP}
The condition Eq.~(\ref{mod-cond})
can be solved for $y$
to yield
\begin{eqnarray}
    y
=
    \tilde{h}_{m+1}(x)
=
      \sum_{i=0}^{m+1} \eps^i \hi(x)
    + {\mathcal O}(\eps^{m+2}) ,
\quad\mbox{for all} \
  x \in K .
\Label{lemma51-exp}
\end{eqnarray}
\end{LEMMA}
Then,
with this first lemma
in hand,
we bootstrap up
from the solution
$y={\tilde h}_{m+1}$
of this modified condition
to find the solution
$y=\hmm$
of the full
$(m+1)-$st
derivative condition,
Eq.~(\ref{cond-alt}).
Specifically,
we show that
their asymptotic expansions
agree up to and including
terms of
${\mathcal O}(\eps^{m+1})$,
\begin{LEMMA}
\label{l-cond}
The condition~(\ref{cond-alt-m+1}),
can be solved for $y$
to yield
\[
    y
=
    \hmm(x)
=
      \sum_{i=0}^{m+1} \eps^i \tilde{h}_{m+1,i}(x)
    + {\mathcal O}(\eps^{m+2}) ,
\quad\mbox{for all} \
  x \in K .
\]
\end{LEMMA}
Given these lemmata
--
the proofs
of which
are given in
appendix~\ref{s-derivcond}
--
Theorem~\ref{t-main1}
follows directly.

\section{Stability analysis of the fixed point
$\hm(x_0)$ \label{s-stage2}}
\begin{figure}[t]
\begin{center}
\leavevmode
\hbox{
  \epsfxsize=6.25in
  \epsfysize=3in
  \epsffile{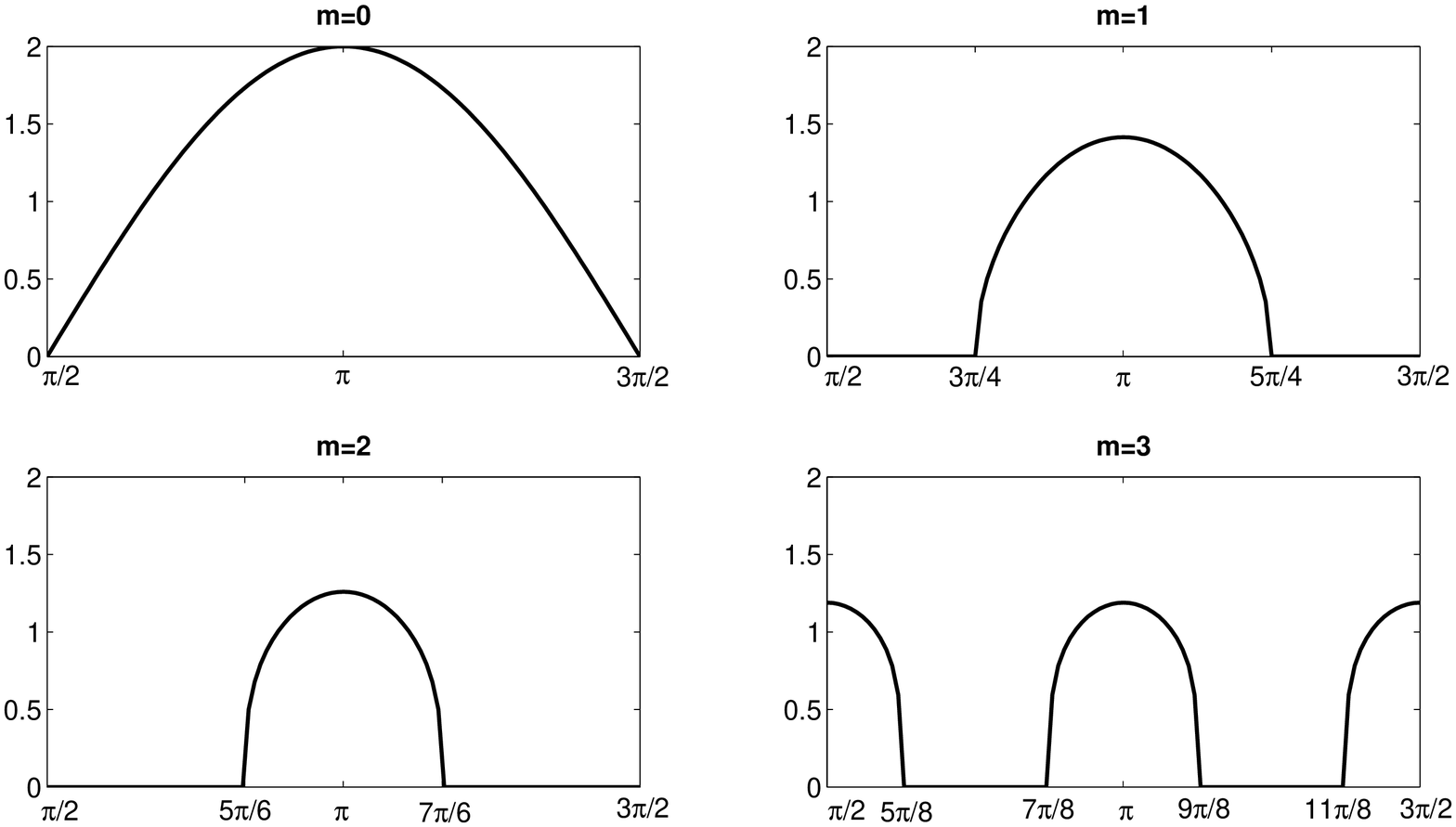}
}
  \caption{\label{f-Hmax-theta}
$H_\ell^{max}$
as a function of
$\theta_\ell \in (\pi/2 ,3\pi/2 )$,
for
$m = 0, 1, 2, 3$.
$H_\ell^{max}$
is measured
in units of
$\eps / \abs{\lambda_\ell}$.
The eigenvalue
$\mu_\ell$
is stable
for all $0 < H < H_\ell^{max}$.}
\end{center}
\end{figure}

\setcounter{equation}{0}
In this section,
we analyze
the stability type
of the fixed point
$y=\hm(x_0)$
of the functional iteration scheme
given by $F_m(y)$.
To fix the notation,
we let
\be
  \sigma(D_yg)_0
=
  \left\{
  \lambda_\ell
  =
  \lambda_{\ell,R}
  +
  i \, \lambda_{\ell,I}
  =
  \vert \lambda_\ell \vert
  {\rm e}^{i \theta_\ell}
  =
  \lambda_{\ell,R}
  (1+i \: {\rm tan}\theta_\ell)
  \, : \,
  \ell = 1 , \ldots , \Nf
  \right\}
\Label{Dyg-spectrum}
\ee
and
remark that
normal attractivity
of the slow manifold
implies that
$\lambda_{\ell,R} < 0$
(equivalently,
$\pi/2 < \theta_\ell < 3\pi/2$)
for all
$\ell = 1 , \ldots , \Nf$.
Then,
we prove the following theorem:
\begin{THEOREM}
\label{t-main2}
For each $m=0,1,\ldots$,
the functional iteration scheme
defined by $F_m$
is stable
if and only if
the following two conditions
are satisfied
for all
$\ell = 1 , \ldots , \Nf$:
\be
  \theta_\ell
\in
  {\mathcal S}_m
\equiv
  \bigcup_{k = 0 , \ldots , m}
\left(
  \frac{2m + 4k + 1}{2(m+1)} \pi
,
  \frac{2m + 4k + 3}{2(m+1)} \pi
\right)
\cap
\left[
\left(
  \frac{\pi}{2} , \frac{3\pi}{2}
\right)
\
{\rm mod} \ 2\pi
\right]
\Label{Sm}
\ee
and
\be
  0
<
  H
<
  \Hmax_\ell
\equiv
  \frac{\eps}{\abs{\lambda_\ell}}
\left[
     2\cos((m+1)(\theta_\ell - \pi))
\right]^{1/(m+1)} .
\Label{Hmax}
\ee
In particular,
if
$\lambda_1 , \ldots , \lambda_\Nf$
are real,
then
the functional iteration is stable
for all $H$ satisfying
\be
  H
<
  \Hmax
\equiv
  2^{1/(m+1)} \ \frac{\eps}{\|D_y g\|_2} .
\Label{H-cond-Re}
\ee
\end{THEOREM}
The graphs of
the stability regions
for $m=0,1,2,3$
are given
in Figure~\ref{f-Hmax-theta}.

We now prove
this theorem.
By definition,
$\hm(x_0)$
is exponentially attracting
if and only if
\be
  \sigma\left(
  \left( DF_m \right)(\hm(x_0))
  \right)
\subset
  {\rm B}(0;1) ,
\Label{DF(m)<B(1,0)}
\ee
where
${\rm B}(0;1)$
denotes the open ball
of radius one
centered at the origin.
To determine the spectrum
of $(DF_m)(\hm(x_0))$,
we use Eq.~(\ref{F(m)-def})
and
Lemma~\ref{l-L(m)}
to obtain
\begin{eqnarray*}
  \left( DF_m \right)(y)
&=&
  I_\Nf
-
  \left( D_y L_m \right)(x_0,y)
\\
&=&
  I_\Nf
-
  \left( -\eps^{-1} H (D_y g)(x_0 , y , 0) \right)^{m+1}
+
  {\mathcal O}\left(\eps, \|g_0(x_0 , y)\| \right) .
\end{eqnarray*}
Letting $y = \hm(x_0)$
in this expression
and
observing that
$\| g_0(x_0 , \hm(x_0)) \|
=
{\mathcal O}(\eps)$
by virtue of
the estimate
$\hm = h_0 + {\mathcal O}(\eps)$
(see Theorem~\ref{t-main1})
and
Eq.~(\ref{inveq-0}),
we obtain
to leading order
\be
  \left( DF_m \right)(\hm(x_0))
=
  I_\Nf
-
  \left( -\eps^{-1} H D_y g \right)_0^{m+1} ,
\Label{DvL(m)}
\ee
where
$z_m = (x_0 , \hm(x_0))$
and
the notation
$(\cdot)_0$
signifies
that the quantity
in parentheses
is evaluated
at the point
$(x_0,\ho(x_0)) \in \Mo$.
Finally, then,
we find
to leading order
\be
  \sigma\left(
  \left( DF_m \right)(\hm(x_0))
  \right)
\!=\!
  \left\{
  \left.
  \mu_\ell
=
  1
-
  \left(
  \abs{\lambda_\ell}\eps^{-1} H
  \right)^{m+1}
  {\rm e}^{i (m+1)(\theta_\ell - \pi)}
  \, \right\vert \,
  \ell = 1 , \ldots , \Nf
  \right\} .
\Label{sigmaDF(m)}
\ee

In view
of Eq.~(\ref{sigmaDF(m)}),
condition~(\ref{DF(m)<B(1,0)})
becomes
\be
  \abs{
  1
-
  \left(
  \abs{\lambda_\ell} \eps^{-1} H
  \right)^{m+1}
  {\rm e}^{i (m+1)(\theta_\ell - \pi)}
  }
\ < \
  1 ,
\quad\mbox{for all} \
  \ell
=
  1 , \ldots , \Nf .
\Label{Cnvrg-cond}
\ee
Here, we note
that higher order terms
omitted from formula~(\ref{sigmaDF(m)})
do not affect stability
for small enough
values of $\eps$,
because
the stability region
${\rm B}(0;1)$ is an open set.
Next,
we study
the circumstances
in which
this stability condition
is satisfied.
This study
naturally splits
into the following
two cases:

\paragraph{Case 1:
The eigenvalues $\lambda_1,\ldots,\lambda_\Nf$ are real.}
This is the case,
for example,
when the fast part
of system~(\ref{z-system})
corresponds to
a spatial discretization
of a self-adjoint operator.
Here,
$\theta_\ell = \pi$
for all
$\ell$,
and thus
condition~(\ref{Cnvrg-cond})
reduces to
\[
  0
\ < \
  \left(
  \abs{\lambda_\ell} \eps^{-1} H
  \right)^{m+1}
\ < \
  2 ,
\quad\mbox{for all} \
  \ell
=
  1 , \ldots , \Nf ,
\]
which further yields
Eq.~(\ref{H-cond-Re}).

\paragraph{Case 2:
Some of the eigenvalues
$\lambda_1,\ldots,\lambda_\Nf$
have nonzero imaginary parts.}
Using Eq.~(\ref{sigmaDF(m)}),
we calculate
\[
  \abs{\mu_\ell}^2
=
  1
+
\left(
  \abs{\lambda_\ell} \eps^{-1} H
\right)^{m+1}
\left[
  \left(
    \abs{\lambda_\ell} \eps^{-1} H
  \right)^{m+1}
-
  2\cos((m+1)(\theta_\ell - \pi))
\right] .
\]
This equation shows that
$\vert\mu_\ell\vert^2$
is
a convex quadratic function
of $H^{m+1}$.
Convexity implies that,
if there exists
a solution $\Hmax_\ell > 0$
to the equation
$\abs{\mu_\ell} = 1$,
then
$\vert\mu_\ell\vert < 1$
for all
$0 < H < \Hmax_\ell$.
Plainly,
$\abs{\mu_\ell} = 1$
implies
\[
\left(
  \abs{\lambda_\ell} \eps^{-1} \Hmax_\ell
\right)^{m+1}
-
  2\cos((m+1)(\theta_\ell - \pi))
=
  0 ,
\]
which yields condition~(\ref{Hmax}).
Further,
the condition that
$\Hmax_1 , \ldots , \Hmax_\Nf$
be real and positive
translates into condition~(\ref{Sm}).
This completes the proof
of Theorem~\ref{t-main2}.

For later comparison
to the results
of numerical simulations,
it is useful to write
formula (\ref{Hmax}) explicitly
for the first several values
of $m$.
For $m = 0$,
formula~(\ref{Hmax})
becomes
\[
   \Hmax_\ell
=
  - \frac{\eps}{\abs{\lambda_\ell}}
    2\cos\theta_\ell ,
\]
see Figure~\ref{f-Hmax-theta}.
We note that
$\Hmax_\ell > 0$
for all
$\theta_\ell \in (\pi/2 , 3\pi/2)$,
and thus
the fixed point
$h_0(x_0)$
is stable
for all
$0 < H < \Hmax$,
where
$\Hmax = \min_\ell(\Hmax_\ell)$.

For $m = 1$,
formula~(\ref{Hmax})
becomes
\[
   \Hmax_\ell
=
  \frac{\eps}{\abs{\lambda_\ell}}
  \sqrt{2\cos(2 \theta_\ell)} ,
\]
see Figure~\ref{f-Hmax-theta}.
We see that,
on
$(\pi/2 , 3\pi/2)$,
$\Hmax_\ell > 0$
\emph{only} if
$\theta_\ell$
lies in the subinterval
$(3\pi/4 , 5\pi/4)$.
Therefore,
the fixed point
$h_1(x_0)$
is stable
if and only if
(i)
$\theta_\ell \in (3\pi/4 , 5\pi/4)$,
for all
$\ell = 1 , \ldots , \Nf$,
and
(ii)
$0 < H < \Hmax = \min_\ell(\Hmax_\ell)$.

For $m = 2$,
formula~(\ref{Hmax})
becomes
\[
   \Hmax_\ell
=
  - \frac{\eps}{\abs{\lambda_\ell}}
    [ 2\cos(3 \theta_\ell) ]^{1/3} ,
\]
see Figure~\ref{f-Hmax-theta}.
Here also,
$\Hmax_\ell > 0$
on
$(\pi/2 , 3\pi/2)$
\emph{only} if
$\theta_\ell$
lies in the subinterval
$(5\pi/6 , 7\pi/6)$.
Thus,
$h_2(x_0)$
is stable
if and only if
(i)
$\theta_\ell \in (5\pi/6 , 7\pi/6)$,
for all
$\ell = 1 , \ldots , \Nf$,
and
(ii)
$0 < H < \Hmax = \min_\ell(\Hmax_\ell)$.

For $m = 3$,
formula~(\ref{Hmax})
becomes
\[
   \Hmax_\ell
=
  \frac{\eps}{\abs{\lambda_\ell}}
  [ 2\cos(4 \theta_\ell) ]^{1/4} ,
\]
see Figure~\ref{f-Hmax-theta}.
We observe that,
on
$(\pi/2 , 3\pi/2)$,
$\Hmax_\ell > 0$
\emph{only} if
$\theta_\ell$
lies in the subdomain
$(\pi/2 , 5\pi/8) \cup (7\pi/8 , 9\pi/8) \cup (11\pi/8 , 3\pi/2)$.
Therefore,
the fixed point
$h_3(x_0)$
is stable
if and only if
(i)
$\theta_\ell
\in
(\pi/2 , 5\pi/8) \cup (7\pi/8 , 9\pi/8) \cup (11\pi/8 , 3\pi/2)$,
for all
$\ell = 1 , \ldots , \Nf$,
and
(ii)
$0 < H < \Hmax = \min_\ell(\Hmax_\ell)$.

\section{Stabilization of the algorithm using RPM \label{s-stage3}}
\setcounter{equation}{0}
In the previous section,
we saw that,
for any $m \ge 1$,
the $m-$th algorithm
in our class of algorithms
may have
a number
of eigenvalues
that
either are unstable
or
have modulus only
slightly less than one.
In this section,
we demonstrate
how
the Recursive Projection Method
(RPM)
of Shroff and Keller
\cite{SK-1993}
may be used
to stabilize
the algorithm
or
to accelerate its convergence
in all such cases.

For the sake of clarity,
we assume that
$(DF_m)(\hm(x_0))$
has $M$ eigenvalues,
labelled
$\{ \mu_1 , \ldots , \mu_M \}$,
that lie outside
the disk
${\rm B}(0 ; 1 - \delta)$,
for some small,
user-specified
$\delta > 0$,
and that
the remaining
$\Nf - M$ eigenvalues
$\{ \mu_{M+1} , \ldots , \mu_\Nf \}$
lie inside it.
We let
$\mathbb{P}$
denote
the maximal invariant subspace
of $(DF_m)(\hm(x_0))$
corresponding to
$\{ \mu_1 , \ldots , \mu_M \}$
and
$P$
denote
the orthogonal projection operator
from $\R^\Nf$
onto that subspace.
Additionally,
we use $\mathbb{Q}$
to denote
the orthogonal complement
of $\mathbb{P}$
in $\R^\Nf$
and
$Q = I_\Nf - P$
to denote
the associated orthogonal projection operator.
These definitions induce
an orthogonal
direct sum decomposition
of $\R^\Nf$,
\[
  \R^\Nf
=
  \mathbb{P}
\oplus
  \mathbb{Q}
=
  P\R^\Nf
\oplus
  Q\R^\Nf ,
\]
and,
as a result,
each
$y \in \R^\Nf$
has a unique decomposition
$y = \tilde{p} + \tilde{q}$,
with
$\tilde{p}
=
Py \in \mathbb{P}$
and
$\tilde{q}
=
Qy\in \mathbb{Q}$.
The fixed point problem
$y = F_m(y)$
may now be written as
\begin{eqnarray}
  \tilde{p}
&=&
  PF_m(\tilde{p} + \tilde{q}) ,
\Label{p-fp}
\\
  \tilde{q}
&=&
  QF_m(\tilde{p} + \tilde{q}) .
\Label{q-fp}
\end{eqnarray}

The fundamental idea
of RPM
is to use
Newton iteration
on Eq.~(\ref{p-fp})
and
functional iteration
on Eq.~(\ref{q-fp}).
In particular,
we decompose
the point
$y^{(1)}$
(which was used
to generate
the sequence
$\{ y^{(r + 1)} \}$
in Eq.~(\ref{seq(m)}))
via
\[
  y^{(1)}
=
  \tilde{p}^{(1)}
+
  \tilde{q}^{(1)}
=
  Py^{(1)}
+
  Qy^{(1)} .
\]
Then,
we apply Newton iteration
on Eq.~(\ref{p-fp})
(starting with $\tilde{p}^{(1)}$)
and
functional iteration
on Eq.~(\ref{q-fp})
(starting with $\tilde{q}^{(1)}$),
\be
\begin{array}{ccl}
  \tilde{p}^{(r + 1)}
&=&
  \tilde{p}^{(r)}
+
  \left[
  I_M
-
  P ( DF_m(\tilde{p}^{(r)} + \tilde{q}^{(r)}) ) P
  \right]^{-1}
  PF_m(\tilde{p}^{(r)} + \tilde{q}^{(r)}) ,
\\
  \tilde{q}^{(r + 1)}
&=&
  QF_m(\tilde{p}^{(r)} + \tilde{q}^{(r)}) .
\end{array}
\Label{RPM}
\ee
The iteration is terminated
when
$\| y^{(r+1)} - y^{(r)} \|
<
\mathrm{TOL}_m$,
for some $r \ge 1$,
as was also the case
with functional iteration.

Application
of Theorem~3.13
from \cite{SK-1993}
directly yields
that the stabilized
(or accelerated)
iterative scheme~(\ref{RPM})
converges
for all initial guesses
$y^{(1)}$
close enough
to the fixed point
$\hm(x_0)$,
as long as
\[
  1
\notin
  \sigma(P ( DF_m(\hm(x_0)) ) P)
=
  \{ \mu_1 , \ldots , \mu_M \} .
\]
In our case,
this condition is satisfied
for all $H > 0$,
because
the fact that
$\Me$
is normally attracting
implies that
each eigenvalue $\lambda_\ell$
of $D_yg$
is bounded away
from zero
uniformly over
the domain $K$
on which
the slow manifold
is defined.
Thus,
the iteration scheme~(\ref{RPM})
converges.

\section{Tuning of the tolerance \label{s-stage4}}
\setcounter{equation}{0}
In this section,
we establish that,
for every
$m=0,1,\ldots$,
$\| y^\#_m - h(x_0) \|
=
{\mathcal O}(\eps^{m+1})$
whenever
$\mathrm{TOL}_m = {\mathcal O} (\eps^{m+1})$.
The value returned
by the functional iteration
is within the tolerance
of the point
on the true slow manifold
for sufficiently small
values of the tolerance.

The brunt
of the analysis
needed to prove
this principal result
involves showing that,
for these small tolerances,
$y^\#_m$
is within the tolerance
of the fixed point,
$\hm(x_0)$.
The desired principal result
is then immediately obtained
by combining this result
with the result of Theorem~\ref{t-main1},
where it was shown that
$\| \hm(x_0) - h(x_0) \|
=
{\mathcal O}(\eps^{m+1})$.

We begin
by observing that
\begin{eqnarray*}
  \| y^\#_m - \hm(x_0) \|
\le
  \| y^\#_m - y^{(r)} \|
+
  \| y^{(r)} - \hm(x_0) \| ,
\quad\mbox{for any} \ r>0 ,
\end{eqnarray*}
by the triangle inequality.
The first term
is ${\mathcal O}(\eps^{m+1})$
by definition,
as long as $r$ is chosen
large enough
so that the stopping criterion,
$\| y^{(r+1)} - y^{(r)} \| < \mathrm{TOL}_m$,
is satisfied.
As to the second term,
we may obtain
the same type
of estimate,
as follows:
First,
\begin{eqnarray*}
  y^{(r+1)} - y^{(r)}
=
  F_m\left( y^{(r)} \right) - y^{(r)}
=
  -L_m\left( x_0 , y^{(r)} \right) ,
\end{eqnarray*}
where we used
Eq.~(\ref{F(m)-def}),
and hence
\[
  L_m\left( x_0 , y^{(r)} \right)
=
  y^{(r)} - y^{(r+1)} .
\]
Second,
$L_m$
is invertible
in a neighborhood
of its fixed point,
by the Implicit Function Theorem,
because
the Jacobian of
$L_m(x_0,\cdot)$
at $\hm(x_0)$ is
\[
  \left( D_y L_m \right)(z_m)
=
  \left( -\eps^{-1} H D_y g \right)_0^{m+1} ,
\]
by Eq.~(\ref{DvL(m)}),
and $\mathrm{det}(D_yg) \ne 0$
since $\Mo$
is normally attracting.
Third,
by combining
these first two observations,
we see that
\[
  y^{(r)}
=
  L^{-1}_m\left( y^{(r)} - y^{(r+1)} \right) ,
\]
where
$L^{-1}_m$ denotes
the local inverse of
$L_m(x_0 , \cdot)$.
Fourth, and finally,
by expanding
$L^{-1}_m$ around zero,
noting that
$L^{-1}_m(0) = \hm(x_0)$,
and
using the triangle inequality,
we obtain
\[
  \| y^{(r)} - \hm(x_0) \|
\le
  \left\| (D_yL^{-1}_m)(0) \right\|
  \left\| y^{(r)} - y^{(r+1)} \right\|
+
  {\mathcal O}\left( \| y^{(r)} - y^{(r+1)} \|^2 \right) .
\]
Recalling the stopping criterion,
we have therefore
obtained the desired bound
on the second term, as well,
\begin{eqnarray*}
  \| y^{(r)} - \hm(x_0) \|
&<&
  \| \left( D_yL^{-1}_m \right)(0) \|
  \mathrm{TOL}_m
+
  {\mathcal O}\left( (\mathrm{TOL}_m)^2 \right) .
\end{eqnarray*}
Hence, the analysis
of this section
is complete.

\section{The effects of differencing \label{s-stage5}}
\setcounter{equation}{0}
In a numerical setting,
the time derivatives of $y$
are approximated,
at each iteration,
by a differencing scheme,
\[
\left(
  \frac{d^{m+1}y}{dt^{m+1}}
\right)(z)
\approx
  \frac{1}{\hat{H}^{m+1}}
  \left( \Delta^{m+1}y \right)(z) ,
\quad\mbox{where}\quad
  z \equiv (x_0 , y)
\quad\mbox{and}\quad
  \hat{H} > 0 .
\]
In this section,
we examine
how the approximation
and convergence results
of Sections~\ref{s-stage1}--\ref{s-stage4}
are affected
by the use
of differencing.
We choose forward differencing,
\be
  \left( \Delta^{m+1}y \right)(z)
=
  \sum_{\ell = 0}^{m+1}
  (-1)^{m + 1 - \ell}
  \left(
  \begin{array}{c}
  m + 1
\\
  \ell
  \end{array}
  \right)
  \phi^y(z ; \ell \hat{H}) ,
\Label{numderiv-def}
\ee
where
$\phi(z ; t)$ is
a (numerically generated) solution
with initial condition $z$,
for concreteness of exposition
and where $\hat{H}$ is a positive,
${\mathcal O}(\eps)$ quantity.
Also,
forward differencing
is directly implementable
in an Equation-Free
or legacy code
setting.

By the Mean Value Theorem,
\begin{eqnarray}
  \left( \Delta^{m+1}y \right)(z)
&=&
  \hat{H}^{m+1}
  \left( \frac{d^{m+1} y}{dt^{m+1}} \right)(z)
+
  \frac{m+1}{2}
  \hat{H}^{m+2}
  \left( \frac{d^{m+2} y}{dt^{m+2}} \right)(\phi(z ; \hat{t}))
\nonumber\\
&=&
  \left(-\frac{1}{\eta}\right)^{m+1}
\left[
  L_m(z)
-
  \frac{m+1}{2\eta}
  L_{m+1}(\phi(z ; \hat{t}))
\right] ,
\Label{numderiv-vs-deriv}
\end{eqnarray}
where
$\eta = H/\hat{H} > 0$
is an ${\mathcal O}(1)$ parameter
available for tuning
and
$\phi(z ; \hat{t})$
is the point
on the solution
$\phi(z ; t)$
at some time
$\hat{t} \in [ 0 , (m+1)\hat{H} ]$.
Thus,
for the $m-$th algorithm,
the approximation of
$d^{m+1}y/dt^{m+1}$
by the above scheme
corresponds
to generating the sequence
$\{ y^{(r)} \vert r = 1, 2, \ldots \}$
using the map
\be
  \hat{F}_m(y)
=
  y
-
  \hat{L}_m(z) ,
\quad
  z
=
  (x_0 , y) ,
\Label{hatF(m)-def-raw}
\ee
where
\be
  \hat{L}_m(z)
=
  (-\eta)^{m+1}
  \left( \Delta^{m+1}y \right)(z)
=
  L_m(z)
-
  \frac{m+1}{2\eta}
  L_{m+1}(\phi(z ; \hat{t})) .
\Label{hatLm}
\ee
Therefore,
by Eq.~(\ref{numderiv-vs-deriv}),
\[
  \hat{F}_m(y)
=
  F_m(y)
+
  \frac{m+1}{2\eta}
  L_{m+1}(\phi(z ; \hat{t})) .
\]

\paragraph{Remark.}
For convenience in the analysis
in this section,
we take
the flow $\phi$
to be
the exact flow
corresponding to Eq.~(\ref{z-system}).
The analysis extends directly
to many problems
for which only
a numerical approximation
of $\phi$
is known.
For example,
if the discretization procedure
admits a smooth error expansion
(such as exists often
for fixed step-size integrators
in legacy codes
or
in the Equation-Free context),
then
the leading order results still hold,
and
the map $\phi$
obtained numerically
is sufficiently accurate
so that the remainder estimates below hold.
In particular,
given a $p$-th order scheme
and
an integration step size $\tilde{h}$,
it suffices to take
$\tilde{h}={\mathcal O}(\eps)$
to guarantee that
the error made in using
the numerically-obtained map $\phi$
is
${\mathcal O}(\eps^p)$.
Of course,
with other integrators,
one could alternatively require
that the timestepper be
${\mathcal O}(\eps^{m+2})$
accurate,
{\it i.e.,} 
of one-higher order of accuracy.

\subsection{Existence of a fixed point $\hat{h}_m(x_0)$ of the map $\hat{F}_m$ \label{ss-stage1}}
In this section,
we establish
that the map
$\hat{F}_m$
has an isolated fixed point
$y = \hat{h}_m(x)$
which differs from $\hm(x_0)$
(and thus also
from $h(x_0)$,
by virtue of Theorem~\ref{t-main1})
only by terms of
${\mathcal O}(\eps^{m+1})$.

The fixed point condition
$\hat{F}_m(x_0 , y) = y$
may be rewritten as
\be
  0
=
  \hat{L}_m(x_0 , y)
=
  L_m(x_0 , y)
-
  \frac{m+1}{2\eta}
  L_{m+1}(\phi(x_0 , y ; \hat{t})) ,
\Label{L-perturb=0}
\ee
where
we combined
Eqs.~(\ref{hatF(m)-def-raw})
and
(\ref{hatLm}).
In order to show
that $\hat{F}_m$
has an isolated fixed point
$\hat{h}_m(x_0)$
which is
${\mathcal O}(\eps^{m+1})-$close
to $\hm(x_0)$,
we need to establish
the validity
of the following
two conditions.
\paragraph{(i)}
The second term
in the right member
of Eq.~(\ref{L-perturb=0})
satisfies
the asymptotic estimate
\be
  \| L_{m+1}(\phi(z_m ; \hat{t})) \|
=
  {\mathcal O}(\eps^{m+1}) ,
\quad\mbox{where} \
  z_m = (x_0 , \hm(x_0)) .
\Label{L(m+1)-cond}
\ee
\paragraph{(ii)}
The Jacobian of
$\hat{L}_m$
satisfies
\be
  {\rm det}
  \left(D_y\hat{L}_m \right)(z_m)
\ne
  0
\quad\mbox{and}\quad
\left\|
  \left(D_y\hat{L}_m \right)(z_m)
\right\|_2
=
  {\mathcal O}(1) .
\Label{DyhatL(m)-cond}
\ee

Let us begin by examining
the term
$L_{m+1}(\phi(z_m ; \hat{t}))$.
Let
$(\hat{x} , \hat{y})
=
\phi(z_m ; \hat{t})$.
Then,
we may write
\begin{eqnarray*}
  L_{m+1}(\phi(z_m ; \hat{t}))
&=&
  L_{m+1}(\hat{x} , \hat{y})
-
  L_{m+1}( \hat{x} , \hmm(\hat{x}) ) ,
\end{eqnarray*}
because
$L_{m+1}(\cdot , \hmm(\cdot) )
\equiv
0$
by the definition
of $L_{m+1}$ and $\hmm$.
Hence,
\be
  \| L_{m+1}(\phi(z_m ; \hat{t})) \|
\le
  \|
  \left( D_yL_{m+1} )(\hat{x} , \hmm(\hat{x} \right) )
\|
\
\|
  \hat{y} - \hmm(\hat{x})
\|
+
\
  {\mathcal O}\left( \| \hat{y} - \hmm(\hat{x}) \|^2 \right) .
\Label{L(m+1)ophi-ineq}
\ee
Now,
$\| (D_yL_{m+1} )(\hat{x} , \hmm(\hat{x})) \|$
is
${\mathcal O}(1)$
by Lemma~\ref{l-L(m)}.
Next,
the triangle inequality yields
\begin{eqnarray*}
  \| \hat{y}  - \hmm(\hat{x}) \|
\le
  \| \hat{y} - h(\hat{x}) \|
+
  \| h(\hat{x}) - \hmm(\hat{x}) \| .
\end{eqnarray*}
The first term
in the right member
remains
${\mathcal O}(\eps^{m+1})$
for all times
$\hat{t} \in [ 0 , (m+1) \hat{H})]$.
Indeed,
the initial condition
$z_m$ is
${\mathcal O}(\eps^{m+1})$-close
to the normally attracting manifold $\Me$.
Thus,
the Fenichel normal form \cite{J-1994}
guarantees that the orbit generated
by this initial condition
remains ${\mathcal O}(\eps^{m+1})$-close
to $\Me$
for ${\mathcal O}(1)$ time intervals.
The second term
in the right member
is also
${\mathcal O}(\eps^{m+1})$,
by Theorem~\ref{t-main1}.
Thus,
$\| \hat{y}  - \hmm(\hat{x}) \|$
is also
${\mathcal O}(\eps^{m+1})$.
Substituting these estimations
into inequality~(\ref{L(m+1)ophi-ineq}),
we obtain that
$\| L_{m+1}(\phi(z_m ; \hat{t})) \|$
is
${\mathcal O}(\eps^{m+1})$
and
condition~(\ref{L(m+1)-cond})
is satisfied.

Next,
we determine
the spectrum of
$(D_y\hat{L}_m)(z_m)$
to leading order
to check
condition~(\ref{DyhatL(m)-cond}).
We will work
with the definition
of $\Delta^{m+1} y$,
Eq.~(\ref{numderiv-def}),
rather than
with formula~(\ref{numderiv-vs-deriv})
which involves
the unknown time $\hat{t}$.
Combining
Eqs.~(\ref{numderiv-def})
and
(\ref{hatF(m)-def-raw}),
we obtain
\[
  \hat{L}_m(z)
=
  \eta^{m+1}
  \sum_{\ell = 0}^{m+1}
  \left(
  \begin{array}{c}
  m + 1
\\
  \ell
  \end{array}
  \right)
  (-1)^\ell
  \phi^y(z ; \ell \hat{H}) .
\]
Differentiating
both members
of this equation
with respect to $y$,
we obtain
\be
  \left (D_y\hat{L}_m \right)(z)
=
  \eta^{m+1}
  \sum_{\ell = 0}^{m+1}
  \left(
  \begin{array}{c}
  m + 1
\\
  \ell
  \end{array}
  \right)
  (-1)^\ell
  (D_y\phi^y)(z ; \ell \hat{H}) .
\Label{DyhatL(m)-aux}
\ee
Next,
$(D_y\phi^y)(z_m ; t)
=
{\rm e}^{(t / \eps) (D_yg)_0}$
to leading order
and
for all
$t$
of ${\mathcal O}(\eps)$
by standard results.
Since
$\ell \hat{H} = {\mathcal O}(\eps)$
for all
$\ell = 0 , 1 , \ldots , (m+1)$,
we may
use this formula
to rewrite
Eq.~(\ref{DyhatL(m)-aux})
to leading order as
\[
  \left (D_y\hat{L}_m \right)(z_m)
=
  \eta^{m+1}
  \sum_{\ell = 0}^{m+1}
  \left(
  \begin{array}{c}
  m + 1
  \\
  \ell
  \end{array}
  \right)
  \left(
  - {\rm e}^{ (\hat{H} / \eps) (D_yg)_0 }
  \right)^\ell
=
  \eta^{m+1}
\left(
  I_\Nf - {\rm e}^{ (\hat{H} / \eps) (D_yg)_0 }
\right)^{m+1} .
\]
Hence,
\be
  \sigma\left( \left(D_y\hat{L}_m \right)(z_m) \right)
=
\left\{
\left.
  \eta^{m+1}
  \left( 1 - {\rm e}^{ \lambda_\ell \hat{H} / \eps } \right)^{m+1}
\right\vert
  \ell = 1 , \ldots , \Nf
\right\} ,
\Label{hatL(m)-spectrum}
\ee
where $z_m = (x_0 , \hm(x_0))$.
This leading order formula
for the elements
of the spectrum
shows that
$(D_y\hat{L}_m)(z_m)$
is ${\mathcal O}(1)$
and non-degenerate
for all
positive ${\mathcal O}(\eps)$ values
of $H$ and $\hat{H}$.
Thus,
condition~(\ref{DyhatL(m)-cond})
is also satisfied.

\subsection{Stability of the fixed point
$\hat h_m(x_0)$ for $\eta = 1$ \label{ss-stage2}}
\begin{figure}[t]
\begin{center}
\leavevmode
\hbox{
  \epsfxsize=6.25in
  \epsfysize=3in
  \epsffile{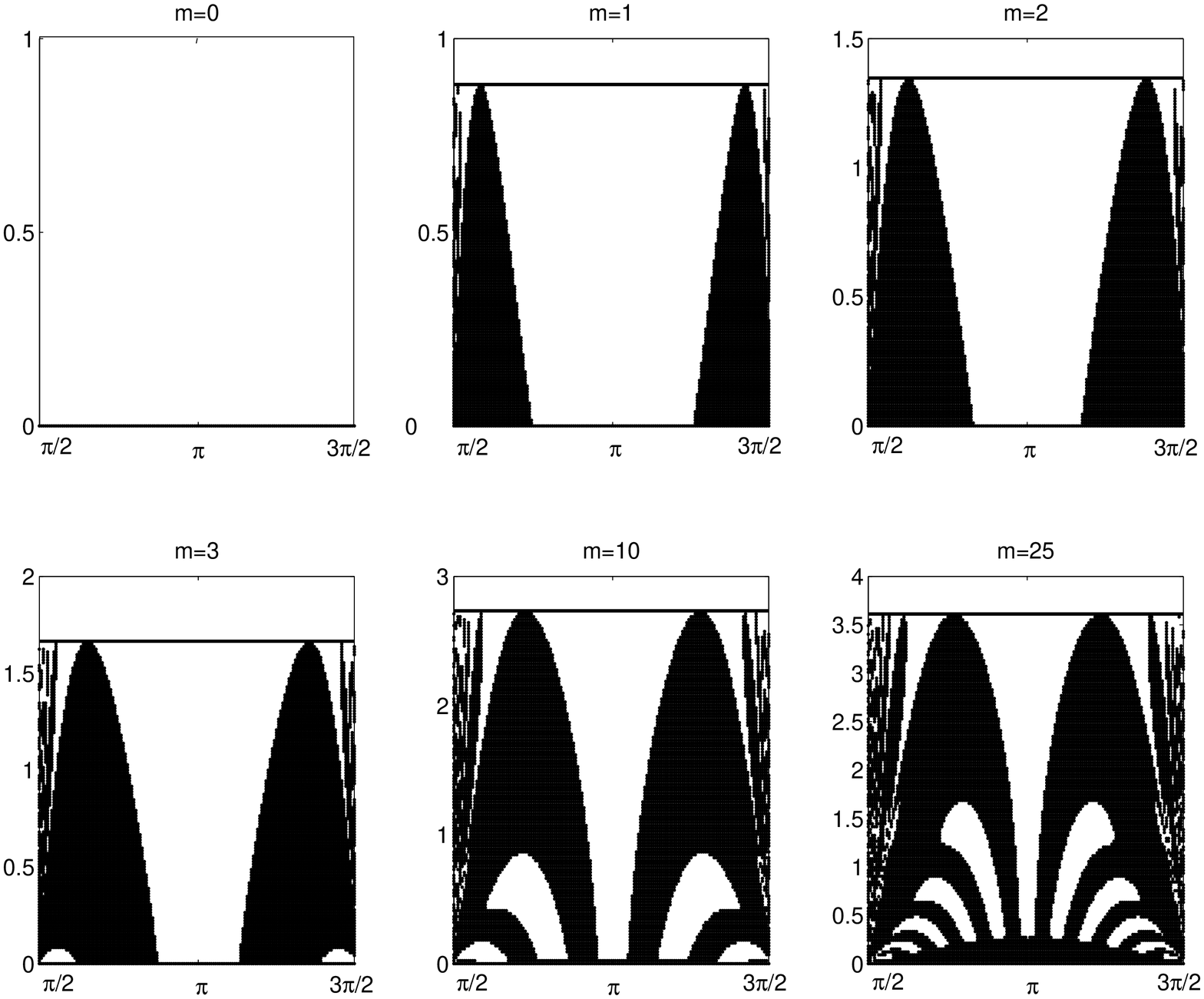}
}
  \caption{\label{f-diff-HvsTheta}
The regions of $H$
for which
$\vert\mu_\ell\vert < 1$
as functions of
$\theta_\ell \in (\pi/2 ,3\pi/2)$.
White corresponds to stability
($\vert\mu_\ell\vert < 1$)
and black to instability
($\vert\mu_\ell\vert > 1$).
$H$ is measured
in units of
$\eps/\abs{\lambda_{\ell,R}}$.
The angle
$\theta_\ell$
takes values
on $(\pi/2 , 3\pi/2)$
and
the black horizontal line
corresponds
to the uniform bound
$H_s(1)$
of Eq.~(\ref{HM}).}
\end{center}
\end{figure}

In this section,
we determine
the stability
of the fixed point
$\hat{h}_m(x_0)$
under functional iteration
using $\hat{F}_m$
in the case that
$\hat{H}=H$.
Our results
for $\hat{H}=H$
are summarized in
the following theorem.
The general case $\hat{H} \ne H$
is treated in the next section,
and
the main result there
is given in
Theorem~\ref{t-main4}.
\begin{THEOREM}
\label{t-main3}
Fix $\eta = 1$.
The functional iteration scheme
defined by $\hat{F}_0$
is unconditionally stable.
For each $m=1,2,\ldots$,
the functional iteration scheme
defined by $\hat{F}_m$
is stable
if and only if,
for each $\ell=1,\ldots,\Nf$,
the pair $(H,\theta_\ell)$
lies in the stability region
the boundary of which
is given by
the implicit equation
\begin{eqnarray}
  1
&=&
  2
  \sum_{j = 1}^{m+1}
  \sum_{k = 1}^{j-1}
  \left(
  \begin{array}{c}
  m + 1
  \\
  j
  \end{array}
  \right)
  \left(
  \begin{array}{c}
  m + 1
  \\
  k
  \end{array}
  \right)
  (-1)^{j+k}
  {\rm e}^{ -(j+k) H_\ell }
  \cos\left( (j-k) H_\ell \tan\theta_\ell \right)
\nonumber\\
&{}&
+
  \sum_{k = 1}^{m+1}
  \left(
  \begin{array}{c}
  m + 1
  \\
  k
  \end{array}
  \right)^2
  {\rm e}^{ -2k H_\ell } ,
\quad\mbox{where}\quad
  H_\ell
=
  - \lambda_{\ell,R} H / \eps> 0 .
\Label{dstabm}
\end{eqnarray}
Here,
the branch of $\arctan$
is chosen so that
$\theta_\ell \in (\pi/2 , 3\pi / 2)$.
In particular,
if
$\lambda_1 , \ldots , \lambda_\Nf$
are real,
then
the functional iteration
is unconditionally stable.
If
at least one
of the eigenvalues
has a nonzero imaginary part,
then
a sufficient and uniform
(in $\theta_1 , \ldots , \theta_\Nf$)
condition for stability
is that
\be
  H
>
  \frac{\eps H_s(1)}{\min_\ell\abs{\lambda_{\ell,R}}} ,
\quad\mbox{where} \
  H_s(1)
=
  - {\rm ln}\left( 2^{1/(m+1)} - 1\right)
\ge
  0 .
\Label{HM}
\ee
\end{THEOREM}
The stability regions
for various values
of $m$
are plotted
in Figure~\ref{f-diff-HvsTheta}.

Following the procedure
used in Section~\ref{s-stage2},
we determine
$\sigma( (D\hat{F}_m )(\hat{h}_m(x_0)) )$
and
examine the circumstances
in which
the stability condition
\be
  \sigma\left(
  \left( D\hat{F}_m \right)(\hat{h}_m(x_0))
  \right)
\subset
  {\rm B}(0;1)
\Label{DhatF(m)<B(1,0)}
\ee
is satisfied.
Equation~(\ref{hatF(m)-def-raw})
yields
\[
  (D\hat{F}_m)(\hat{h}_m(x_0))
=
  I_\Nf
-
  (D_y\hat{L}_m)(x_0 , \hat{h}_m(x_0))
\]
and thus also
\[
  \{ \hat{\mu}_\ell \}
\equiv
  \sigma\left( \left( D_y\hat{F}_m \right)\left( \hat{h}_m(x_0) \right) \right)
=
  1
-
  \sigma\left( \left( D_y\hat{L}_m \right)\left( x_0 , \hat{h}_m(x_0) \right) \right) .
\]
Since
$\hat{h}_m(x_0)$
differs from
$h_m(x_0)$
only at terms of
${\mathcal O}(\eps^{m+1})$,
$(D_y\hat{L}_m)(x_0 , \hat{h}_m(x_0))$
also differs from
$(D_y\hat{L}_m)(x_0 , \hm(x_0))$
only
at terms of
${\mathcal O}(\eps^{m+1})$.
Thus,
Eq.~(\ref{hatL(m)-spectrum})
yields,
to leading order
and
for
$\ell = 1 , \ldots , \Nf$,
\be
  \hat{\mu}_\ell
=
  1 - \left( 1 - {\rm e}^{ \lambda_\ell H / \eps } \right)^{m+1}
=
  \sum_{k = 1}^{m+1}
  \left(
  \begin{array}{c}
  m + 1
  \\
  k
  \end{array}
  \right)
  (-1)^{k+1}
  {\rm e}^{ k \lambda_\ell H / \eps } .
\Label{hatF(m)-spectrum}
\ee
Recalling Eq.~(\ref{Dyg-spectrum})
and
defining
$H_\ell
=
- \lambda_{\ell,R} H / \eps$,
we rewrite
Eq.~(\ref{hatF(m)-spectrum})
in the form
\be
  \hat{\mu}_\ell
=
  \sum_{k = 1}^{m+1}
  \left(
  \begin{array}{c}
  m + 1
  \\
  k
  \end{array}
  \right)
  (-1)^{k+1}
  {\rm e}^{ -k H_\ell (1 + i \tan\theta_\ell) } .
\Label{hatmu}
\ee
The stability condition~(\ref{DhatF(m)<B(1,0)})
becomes, then,
\be
  \abs{
  \hat{\mu}_\ell
  }
=
  \abs{
  \sum_{k = 1}^{m+1}
  \left(
  \begin{array}{c}
  m + 1
  \\
  k
  \end{array}
  \right)
  (-1)^{k+1}
  {\rm e}^{ -k H_\ell (1 + i \tan\theta_\ell) }
  }
\ < \
  1 ,
\quad\mbox{for all} \
  \ell
=
  1 , \ldots , \Nf .
\Label{Cnvrg-cond-hat}
\ee
As in Section~\ref{s-stage2},
we distinguish
two cases.

\paragraph{Case 1:
All
of the eigenvalues
of $(D_y g)_0$
are real.}
Then,
$\theta_\ell = \pi$
for all
$\ell = 1 , \ldots , \Nf$,
and
hence
Eq.~(\ref{hatmu})
becomes
\[
  \hat{\mu}_\ell
=
  \sum_{k = 1}^{m+1}
  \left(
  \begin{array}{c}
  m + 1
  \\
  k
  \end{array}
  \right)
  (-1)^{k+1}
  {\rm e}^{ -k H_\ell }
=
  1 - (1 - {\rm e}^{ -H_\ell })^{m+1} .
\]
Thus,
the spectrum of
$( D_y\hat{F}_m )(\hat{h}_m(x_0))$
is contained in
$(0 , 1)$
for all
positive
${\mathcal O}(\eps)$
values of $H$.
Equivalently,
the fixed point $\hat{h}_m(x_0)$
is
\emph{unconditionally stable}
for these values of $H$.

These results may be interpreted
both in the context
of the $m$-th iterative algorithm
for each fixed $m$,
as well as
in the context
of using the algorithms
as an integrated class.
In particular,
for each fixed $m$,
the rate of convergence
to the fixed point
of the $m$-th algorithm
increases as $H$ increases.
Also, for any fixed
iterative step size $H$,
the rate of convergence
of the $m$-th algorithm
to its fixed point
decreases as the order, $m$,
of the iterative algorithm increases. 
This information is important
for determining how large an $m$ one should use,
especially when using the algorithms
as an integrated class.

\begin{figure}[t]
\begin{center}
\leavevmode
\hbox{
  \epsfxsize=6.25in
  \epsfysize=3in
  \epsffile{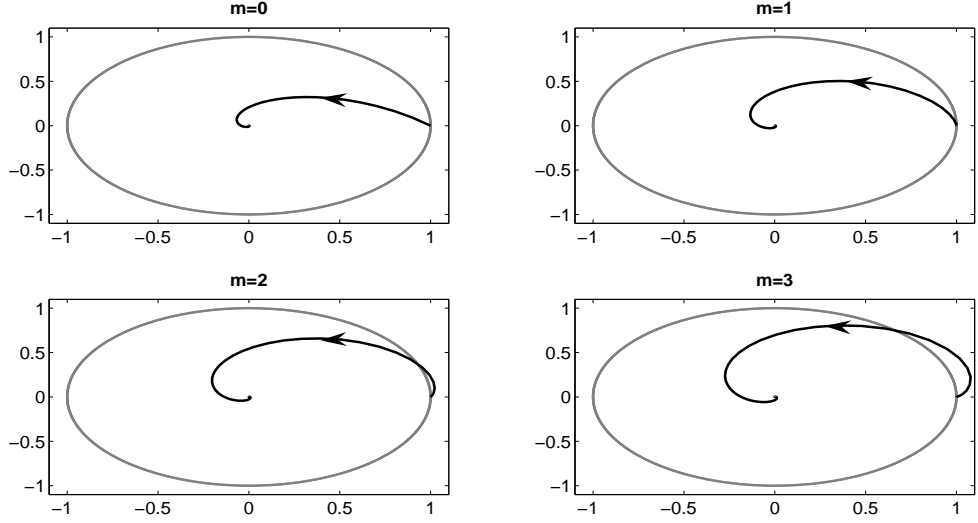}
}
  \caption{\label{f-diff-evalue}
The eigenvalue $\hat{\mu}_\ell$
for values of $H$
between zero and $100\eps$.
The thick line
denotes the boundary
of the stability region
(\emph{i.e.},
the unit circle).
The eigenvalue
$\lambda_\ell$
was taken to be
$-1 + i$
for each one
of the graphs.
The arrow points
to increasing values of $H$.}
\end{center}
\end{figure}

\paragraph{Case 2:
Some of the eigenvalues
of $(D_yg)_0$
have nonzero imaginary parts.}
When this is the case,
some of the eigenvalues
may be unstable
for certain values of $H$.
Figure~\ref{f-diff-evalue} demonstrates this:
in it,
we have drawn
the complex eigenvalue
$\hat{\mu}_\ell$
for various values of $H$
and
for $m = 0,1,2,3$.
Plainly,
$\hat{\mu}_\ell$
is unstable
for $m > 0$
and
for $H$ small enough,
as $\abs{\hat{\mu_\ell}} > 1$.
We determine
the stability regions
in the
$(\theta_\ell , H_\ell)-$plane
as functions of $m$.

First,
we derive
the uniform bound
(\ref{HM}).
Using formula~(\ref{hatmu}),
we calculate
\be
  \abs{\hat{\mu}_\ell}
\le
  \sum_{k = 1}^{m+1}
  \left(
  \begin{array}{c}
  m + 1
  \\
  k
  \end{array}
  \right)
  {\rm e}^{ -k H_\ell }
=
  ( 1 + {\rm e}^{-H_\ell} )^{m+1}
-
  1 ,
\Label{hatmu-mod-estim}
\ee
and thus
$\abs{\hat{\mu}_\ell} < 1$,
for all
$H_\ell > H_s(1)$.
Recalling that
$H_\ell
=
- \lambda_{\ell,R} H / \eps$,
we conclude that
all of the eigenvalues
$\hat{\mu}_\ell$
lie in the unit disk
(equivalently,
the $m-$th algorithm
is stable)
for all
${\mathcal O}(\eps)$ values
of $H$
greater than
$\eps H_s(1) / \min_\ell\abs{\lambda_{\ell,R}}$,
irrespective
of the values of
$\theta_1 , \ldots , \theta_\Nf$.
This is demonstrated
in Figure~\ref{f-diff-HvsTheta}.

Next,
we derive formulae
which describe exactly
the stability regions.
For $m = 0$,
Eq.~(\ref{HM}) yields
$H_s(1) = 0$.
Thus,
$\abs{\hat{\mu}_\ell} < 1$
for all  positive
${\mathcal O}(\eps)$
values of $H$
and
for all
$\ell = 1 , \ldots , \Nf$.
As a result,
the fixed point $\hat{h}_0(x_0)$
is
\emph{unconditionally stable}
for positive,
${\mathcal O}(\eps)$
values of $H$,
see also
Figure~\ref{f-diff-HvsTheta}.

\begin{figure}[t]
\begin{center}
\leavevmode
\hbox{
  \epsfxsize=6.25in
  \epsfysize=3in
  \epsffile{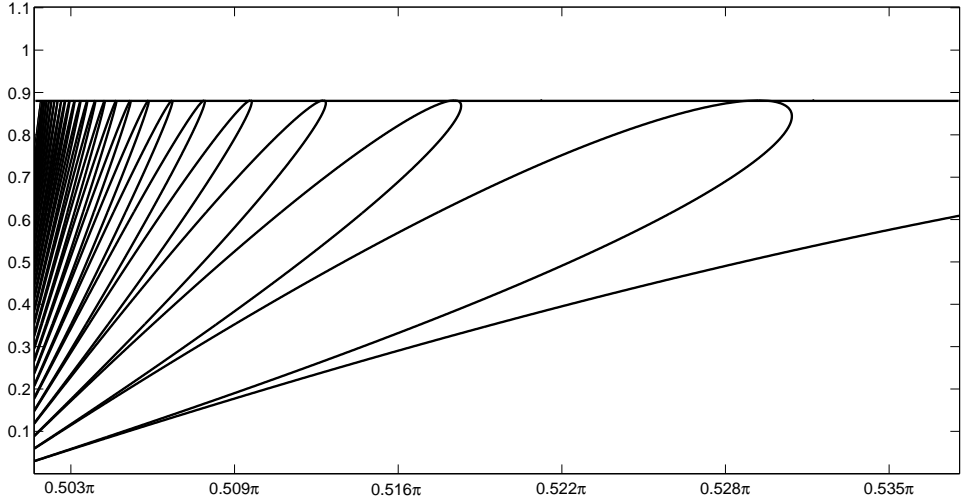}
}
  \caption{\label{f-finestructure}
The fine structure
of the stability region
depicted in Figure~\ref{f-diff-HvsTheta}
(with $m = 1$)
close to $\pi/2$.
The exterior
of the lobes
is part of
the stability region.
}
\end{center}
\end{figure}

For $m = 1$,
Eq.~(\ref{hatmu})
becomes
\[
  \hat{\mu}_\ell
=
  2 {\rm e}^{ -H_\ell (1 + i \tan\theta_\ell) }
-
  {\rm e}^{ -2 H_\ell (1 + i \tan\theta_\ell) } .
\]
Writing
$\overline{\hat{\mu}_\ell}$
for
the complex conjugate
of $\hat{\mu}_\ell$,
then,
we calculate
\be
  \abs{\hat{\mu}_\ell}^2
=
  \hat{\mu}_\ell
\
  \overline{\hat{\mu}_\ell}
=
  4 {\rm e}^{ -2 H_\ell }
-
  4 {\rm e}^{ -3 H_\ell } \cos(H_\ell \, \tan\theta_\ell)
+
  {\rm e}^{ -4 H_\ell } .
\ee
Using this formula,
we recast
the stability condition
(\ref{Cnvrg-cond-hat})
into the form
\[
  4 {\rm e}^{ -2 H_\ell }
-
  4 {\rm e}^{ -3 H_\ell }
  \cos(H_\ell \tan\theta_\ell)
+
  {\rm e}^{ -4 H_\ell }
\ < \
  1 .
\]
In particular,
the boundary
of the stability region
can be obtained
by equating
the expression
in the left member
of this inequality
to one
and
solving for $\theta_\ell$,
to obtain
\[
  \theta_\ell
=
  \arctan
\left(
  H_\ell^{-1}
\left[
  \arccos
\left[
  \frac{1}{4}
  {\rm e}^{ -H_\ell }
+
  {\rm e}^{ H_\ell }
-
  \frac{1}{4}
  {\rm e}^{ 3 H_\ell }
\right]
+
  2 k \pi
\right]
\right) .
\]
Here,
$k \in {\mathbf Z}$
and
the branch of
$\arctan$
is chosen
so that
$\theta_\ell \in (\pi/2 , 3\pi / 2)$.
We have plotted
the stability region
in Figure~\ref{f-diff-HvsTheta}.
We also note here
that
the boundary
of the stability region
close to
$\pi/2$ and to $3\pi/2$
has fine structure,
see Figure~\ref{f-finestructure}.

For a general
value of $m$,
the stability condition
(\ref{Cnvrg-cond-hat}) is
\[
   \abs{\hat{\mu}_\ell}
=
  \abs{
  \sum_{k = 1}^{m+1}
  \left(
  \begin{array}{c}
  m + 1
  \\
  k
  \end{array}
  \right)
  (-1)^{k+1}
  {\rm e}^{ -k H_\ell (1 + i \tan\theta_\ell) }
  }
\ < \
  1 ,
\quad\mbox{for all} \
  \ell
=
  1 , \ldots , \Nf .
\]
Now,
using
Eq.~(\ref{hatmu}),
we calculate
\begin{eqnarray*}
  \abs{\hat{\mu}_\ell}^2
&=&
  \hat{\mu}_\ell
\
  \overline{\hat{\mu}_\ell}
\nonumber\\
&=&
  \sum_{j = 1}^{m+1}
  \sum_{k = 1}^{m+1}
  \left(
  \begin{array}{c}
  m + 1
  \\
  j
  \end{array}
  \right)
  \left(
  \begin{array}{c}
  m + 1
  \\
  k
  \end{array}
  \right)
  (-1)^{j+k}
  {\rm e}^{ -(j+k) H_\ell }
  {\rm e}^{ i (j-k) H_\ell \tan\theta_\ell }
\nonumber\\
&=&
  2
  \sum_{j = 1}^{m+1}
  \sum_{k = 1}^{j-1}
  \left(
  \begin{array}{c}
  m + 1
  \\
  j
  \end{array}
  \right)
  \left(
  \begin{array}{c}
  m + 1
  \\
  k
  \end{array}
  \right)
  (-1)^{j+k}
  {\rm e}^{ -(j+k) H_\ell }
  \cos\left( (j-k) H_\ell \tan\theta_\ell \right)
\nonumber\\
&{}&
+
  \sum_{k = 1}^{m+1}
  \left(
  \begin{array}{c}
  m + 1
  \\
  k
  \end{array}
  \right)^2
  {\rm e}^{ -2k H_\ell } .
\end{eqnarray*}
Equation~(\ref{dstabm})
now follows directly.

\subsection{Stability of the fixed point
$\hat h_m(x_0)$ for $\eta \ne 1$ \label{ss-stage3}}
\begin{figure}[t]
\begin{center}
\leavevmode
\hbox{
  \epsfxsize=6.25in
  \epsfysize=3in
  \epsffile{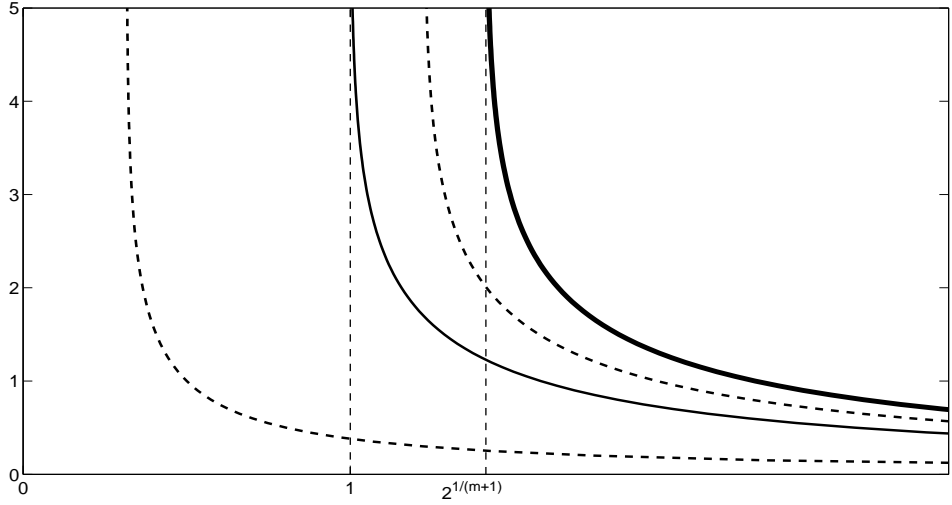}
}
  \caption{\label{f-avshatH}
The stability region
in the $(\eta,\hat{H}_\ell)-$plane
together with the level curves
$\hat{\mu}_\ell(\eta , \hat{H}_\ell) = -1$
(thick curve),
$\hat{\mu}_\ell(\eta , \hat{H}_\ell) = 0$
(solid curve in the middle),
$\hat{\mu}_\ell(\eta , \hat{H}_\ell) = 1$
(union of the two semiaxes).
The dashed level curves
to the right and left
of the level curve
$\hat{\mu}_\ell = 0$
correspond to representative
positive and negative
values of $\hat{\mu}_\ell$,
respectively.
The eigenvalue
$\hat{\mu}_\ell$
is stable
for all pairs $(\eta,\hat{H}_\ell)$
to the left of the level curve
$\hat{\mu}_\ell = -1$.}
\end{center}
\end{figure}

\begin{figure}[t]
\begin{center}
\leavevmode
\hbox{
  \epsfxsize=6.25in
  \epsfysize=3in
  \epsffile{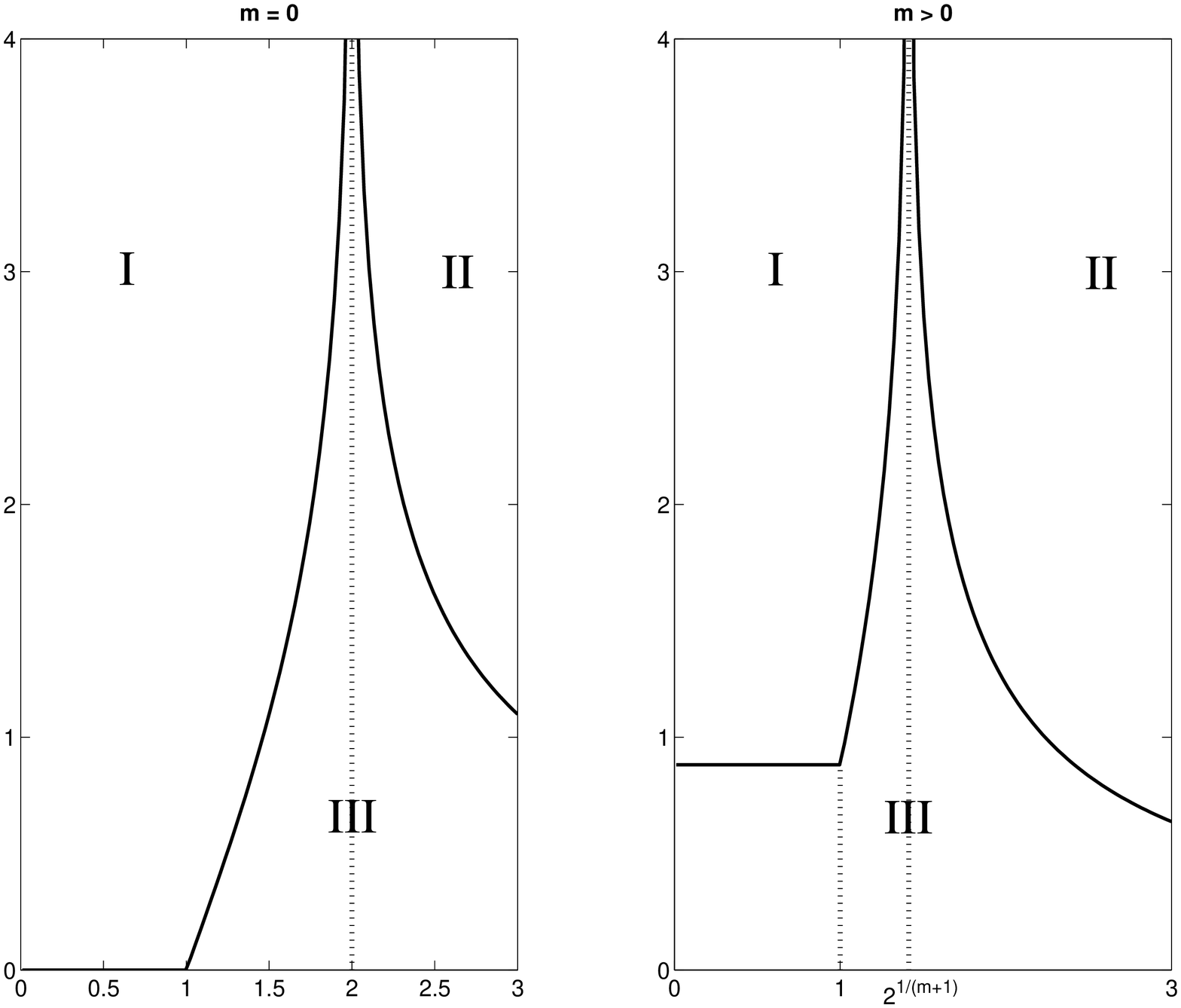}
}
  \caption{\label{f-UnifBd}
The stability regions
in the $(\eta,\hat{H}_\ell)-$plane
for
$m=0$ (left panel)
and
$m = 1 , 2 , \ldots$ (right panel).
The eigenvalue
$\hat{\mu}_\ell$
is stable
in region I,
unstable
in region II,
and
its stability type
is $\theta_\ell-$dependent
in region III.}
\end{center}
\end{figure}

In this section,
we determine
the stability
of the fixed point
$\hat{h}_m(x_0)$
for $\hat{H} \ne H$.
We define the function
\be
  \hat{H}_m(\eta)
=
\left\{
\begin{array}{lcr}
  - {\rm ln}
\left(
  2^{1/(m+1)} - 1
\right) ,
&
\quad
\mbox{if}
\quad
&
  0 < \eta \le 1 ,
\\
  - {\rm ln}
\abs{
  2^{1/(m+1)}/\eta - 1
} ,
&
\quad
\mbox{if}
\quad
&
  \eta > 1 .
\end{array}
\right.
\Label{Hm}
\ee
Our results are summarized
in the following theorem.
\begin{THEOREM}
\label{t-main4}
Fix $\eta > 0$.
For each $m=0,1,2,\ldots$,
the functional iteration scheme
defined by $\hat{F}_m$
is stable
if and only if,
for each $\ell=1,\ldots,\Nf$,
the pair $(\hat{H},\theta_\ell)$
lies in the stability region
the boundary of which
is given by
the implicit equation
\begin{eqnarray}
  1
&=&
  2 \eta^{2(m+1)}
  \sum_{j = 1}^{m+1}
  \sum_{k = 1}^{j-1}
  \left(
  \begin{array}{c}
  m + 1
  \\
  j
  \end{array}
  \right)
  \left(
  \begin{array}{c}
  m + 1
  \\
  k
  \end{array}
  \right)
  (-1)^{j+k}
  {\rm e}^{ -(j+k) \hat{H}_\ell }
  \cos\left( (j-k) \hat{H}_\ell \tan\theta_\ell \right)
\nonumber\\
&{}&
+
  2 \eta^{m+1} \left( \eta^{m+1} - 1 \right)
  \sum_{k = 1}^{m+1}
  \left(
  \begin{array}{c}
  m + 1
  \\
  k
  \end{array}
  \right)
  (-1)^k
  {\rm e}^{ -k \hat{H}_\ell }
  \cos\left(k \hat{H}_\ell \tan\theta_\ell \right)
\nonumber\\
&{}&
+
  \eta^{2(m+1)}
  \sum_{k = 1}^{m+1}
  \left(
  \begin{array}{c}
  m + 1
  \\
  k
  \end{array}
  \right)^2
  {\rm e}^{ -2k \hat{H}_\ell }
+
  \left( \eta^{m+1} - 1 \right)^2 ,
\Label{dstabm-a}
\end{eqnarray}
where
$\hat{H}_\ell
=
- \lambda_{\ell,R} \hat{H} / \eps> 0$.
Here,
the branch of $\arctan$
is chosen so that
$\theta_\ell \in (\pi/2 , 3\pi / 2)$.
In particular:
\\
{\rm (i)}
Assume that
${\rm Im}(\lambda_\ell)
=
0$,
for all
$\ell = 1 , \ldots , \Nf$.
If
$0 < \eta < 2^{1/(m+1)}$,
then
the functional iteration
is unconditionally stable.
If
$\eta > 2^{1/(m+1)}$,
then
the functional iteration
is stable
if and only if
\be
  0
  <
  \hat{H}
<
\frac{
  \eps \hat{H}_m(\eta)
}{
  \max_\ell\abs{\lambda_{\ell,R}}
} .
\Label{H<Hm}
\ee
\\ {\rm (ii)}
Assume that
at least one of
${\rm Im}(\lambda_1),
\ldots,
{\rm Im}(\lambda_\Nf)$
is nonzero.
If
$0 < \eta < 2^{1/(m+1)}$,
then
a sufficient and uniform
(in $\theta_1 , \ldots , \theta_\Nf$)
condition for stability
is
\be
  \hat{H}
>
  \frac{\eps \hat{H}_m(\eta)}
       {\min_\ell\abs{\lambda_{\ell,R}}} .
\Label{H>Hm-stab}
\ee
If
$\eta > 2^{1/(m+1)}$,
the functional iteration
is unstable
for any
$\theta_1 , \ldots , \theta_\Nf$
and
for all
\be
  \hat{H}
>
  \frac{\eps \hat{H}_m(\eta)}
       {\max_\ell\abs{\lambda_{\ell,R}}} .
\Label{H>Hm-instab}
\ee
\end{THEOREM}
These results
are demonstrated
in Figures~\ref{f-avshatH}
and \ref{f-UnifBd}.

As in Section~\ref{ss-stage2},
we determine when
the stability condition (\ref{DhatF(m)<B(1,0)})
holds.
The analogue of
Eqs.~(\ref{hatF(m)-spectrum})
and
(\ref{hatmu})
in this case is,
to leading order
and
for
$\ell = 1 , \ldots , \Nf$,
\begin{eqnarray}
  \hat{\mu}_\ell
&=&
  1
-
  \eta^{m+1}
\left(
  1 - {\rm e}^{ \lambda_\ell \hat{H} / \eps }
\right)^{m+1}
=
  1
-
  \eta^{m+1}
\left(
  1 - {\rm e}^{ -\hat{H}_\ell (1 + i \tan\theta_\ell) }
\right)^{m+1} .
\Label{hatF(m)-spectrum-a}
\end{eqnarray}
The stability condition~(\ref{DhatF(m)<B(1,0)})
becomes, then,
\be
\abs{
  \hat{\mu}_\ell
}
=
\abs{
  1
-
  \eta^{m+1}
\left(
  1 - {\rm e}^{ -\hat{H}_\ell (1 + i \tan\theta_\ell) }
\right)^{m+1}
}
\ < \
  1 ,
\quad\mbox{for all}\quad
  \ell
=
  1 , \ldots , \Nf .
\Label{Cnvrg-cond-hat-a}
\ee
Here also,
we distinguish
two cases.

\paragraph{Case 1:
All
of the eigenvalues
of $(D_y g)_0$
are real.}
Then,
$\theta_\ell = \pi$
for all
$\ell = 1 , \ldots , \Nf$,
and
hence
Eq.~(\ref{Cnvrg-cond-hat-a})
becomes
\[
  \hat{\mu}_\ell
=
  1
-
  \eta^{m+1}
  (1 - {\rm e}^{ -\hat{H}_\ell })^{m+1} .
\]
Plainly,
the condition
$\hat{\mu}_\ell < 1$
is satisfied
for all positive
$\hat{H}_\ell$ and $\eta$.
Next,
solving this equation
for $\eta$,
we obtain an equation
for the level curve
$\hat{\mu}_\ell = constant$,
\[
  \eta
=
  \frac{
\left(
  1
-
  \hat{\mu}_\ell
\right)^{1/(m+1)}}{
  1
-
  {\rm e}^{ -\hat{H}_\ell }} .
\]
For $0 < \eta < 2^{1/(m+1)}$
and
for all
${\mathcal O}(\eps)$
and
positive values of $\hat{H}$,
we obtain
$\hat{\mu}_\ell > -1$
(and thus
the eigenvalue $\hat{\mu}_\ell$
is stable),
see Fig.~\ref{f-avshatH}.
Therefore,
$\sigma(( D_y\hat{F}_m )(\hat{h}_m(x_0)))
\subset
(-1 , 1)$,
and
the fixed point $\hat{h}_m(x_0)$
is
\emph{unconditionally stable}.

For $\eta > 2^{1/(m+1)}$,
we obtain the condition
$0 < \hat{H}_\ell < \hat{H}_m(\eta)$,
and
Eq.~(\ref{H<Hm})
follows directly.
Finally,
we note that,
for a fixed value of $\eta$
and
as $\hat{H} \to \infty$,
the spectrum clusters around
$1 - \eta^{m+1}$.
Thus,
the choice $\eta = 1$
is optimal
in the sense that
large values of $\hat{H}$
bring the spectrum
closer to zero.

\paragraph{Case 2:
Some of the eigenvalues
of $(D_yg)_0$
have nonzero imaginary parts.}
In this case,
some of the eigenvalues
may become unstable
for certain combinations
of $\eta$ and $\hat{H}$,
as
our analysis in
Section~\ref{ss-stage2}
also showed.

First,
we consider the case
$0 < \eta < 2^{1/(m+1)}$
and
derive
the uniform bound
(\ref{H>Hm-stab}).
Using formula~(\ref{hatF(m)-spectrum-a})
and
working as in Eq.~(\ref{hatmu-mod-estim}),
we estimate
\[
  \abs{\hat{\mu}_\ell}
\le
\abs{
  1 - \eta^{m+1}
}
+
  \eta^{m+1}
\left[
  ( 1 + {\rm e}^{-\hat{H}_\ell} )^{m+1}
-
  1
\right] .
\]
Hence
\[
  \abs{\hat{\mu}_\ell}
\le
\left\{
\begin{array}{ll}
  1
+
  \eta^{m+1}
\left[
  ( 1 + {\rm e}^{-\hat{H}_\ell} )^{m+1} - 2
\right] ,
&
\quad\mbox{for}\quad
  0 < \eta \le 1 ,
\\
  \eta^{m+1}
  ( 1 + {\rm e}^{-\hat{H}_\ell} )^{m+1}
-
  1 ,
&
\quad\mbox{for}\quad
  \eta > 1 .
\end{array}
\right.
\]
Combining these inequalities
with the stability condition
$\abs{\hat{\mu}_\ell} < 1$,
we obtain
the sufficient condition
$\hat{H}_\ell > \hat{H}_m(\eta)$,
where
$\hat{H}_m(\eta)$
is the uniform bound~(\ref{Hm})
(see also Fig.~\ref{f-UnifBd}).
Recalling that
$\hat{H}_\ell
=
- \lambda_{\ell,R} \hat{H} / \eps$,
we conclude that,
if
condition (\ref{H>Hm-stab})
is satisfied,
then
$\sigma(( D_y\hat{F}_m )(\hat{h}_m(x_0)))
\subset
{\rm B}(0;1)$,
and hence
the $m-$th algorithm
is stable.

Next,
we consider the case
$\eta > 2^{1/(m+1)}$
and
derive the uniform bound
(\ref{H>Hm-instab}).
Equation~(\ref{hatF(m)-spectrum-a})
yields
\[
  \abs{1 - \hat{\mu}_\ell}
\ge
  \eta^{m+1}
\left(
  1 - \abs{{\rm e}^{ - \hat{H}_\ell }{\rm e}^{ i \hat{H}_\ell \tan\theta_\ell }}
\right)^{m+1}
\ge
  \eta^{m+1}
\left(
  1 - {\rm e}^{ - \hat{H}_\ell }
\right)^{m+1} .
\]
Thus,
$\abs{1 - \hat{\mu}_\ell} > 2 $,
for
$\eta > 2^{1/(m+1)}$
and
$\hat{H}_\ell > \hat{H}_m(\eta)$,
and therefore
\[
  \abs{\hat{\mu}_\ell}
\ge
  \abs{\abs{1 - \hat{\mu}_\ell} - 1}
>
  1 ,
\]
Hence,
$\hat{\mu}_\ell$ is unstable.
\begin{figure}[t]
\begin{center}
\leavevmode
\hbox{
  \epsfxsize=6.25in
  \epsfysize=3in
  \epsffile{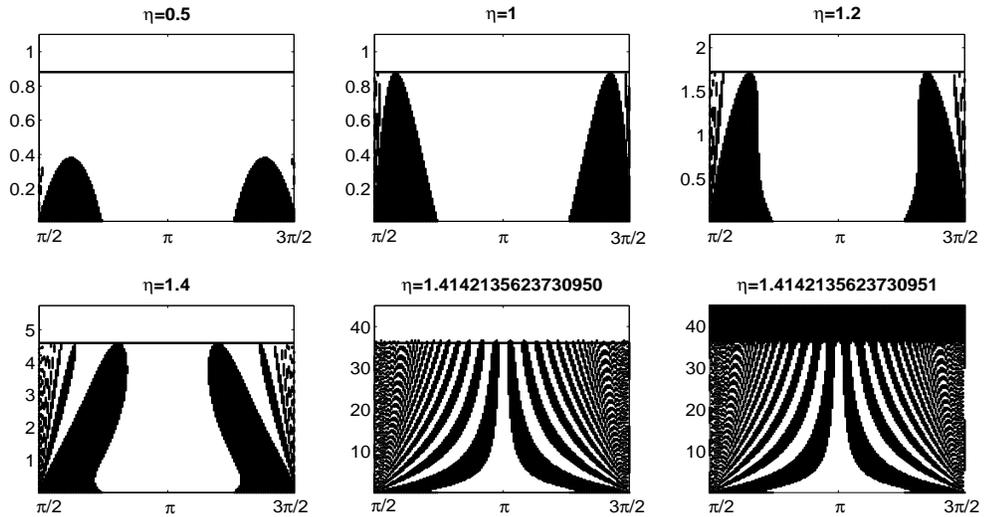}
}
  \caption{\label{f-Bif}
The stability region
in the $(\eta,\hat{H}_\ell)-$plane
for $m=1$
and for various values
of $\eta$.
The last two values for $\eta$
are just below
and just above
the value
$2^{1/(m+1)}=\sqrt{2}$.
White denotes stability
and
black denotes instability.}
\end{center}
\end{figure}

\paragraph{Remark.}
Conditions (\ref{H>Hm-stab})
and
(\ref{H>Hm-instab})
may be interpreted
by means of
the fact that
$\sigma(( D_y\hat{F}_m )(\hat{h}_m(x_0)))$
clusters around
$1 - \eta^{m+1}$
as $\hat{H}\to\infty$.
For
$0 < \eta < 2^{1/(m+1)}$,
there holds that
$-1 < 1 - \eta^{m+1} < 1$.
Thus,
for $\hat{H}$ large enough,
the eigenvalues are contained
in the unit disk.
On the contrary,
$1 - \eta^{m+1} < -1$
for
$\eta > 2^{1/(m+1)}$,
and thus
the eigenvalues lie outside
the unit disk
for $\hat{H}$ large enough.

Finally,
formula~(\ref{dstabm-a})
describing the stability region
may be derived
in a manner
entirely analogous
to that used
to derive Eq.~(\ref{dstabm}).


\section{Conclusions and Discussion}
\setcounter{equation}{0}
In this article,
we characterized
the accuracy and convergence properties
of the class of iterative algorithms
introduced in \cite{GKKZ-2005}
for explicit fast-slow systems (\ref{z-system}).
The $m$-th member of the class
corresponds to
a functional iteration scheme
to solve
the $(m+1)-$st derivative condition
(\ref{deriv-cond}).
We showed that
this condition has
an isolated solution,
which corresponds to
a fixed point
of this $m$-th member
and
which is accurate
up to and including
terms of
${\cal O}(\eps^{m})$,
see Theorem~\ref{t-main1}.
Also,
we derived explicit formulae
for the domain of convergence
of the functional iteration,
both
in the case where
analytical formulae
for the $(m+1)-$st derivative
are used
(see Theorem~\ref{t-main2})
and
in the case where
the $(m+1)-$st derivatives
are estimated through
a forward difference scheme
(see Theorem~\ref{t-main3}).
These convergence results
are illustrated in Figures~\ref{f-Hmax-theta},
\ref{f-diff-HvsTheta},
and \ref{f-finestructure}.
Further,
we demonstrated how
the Recursive Projection Method
may be used to stabilize
the functional iteration
in all cases when it is unstable
or
to accelerate its convergence
in those cases where the convergence is slow.

An extension
of the analysis presented here
to more general multiscale systems (\ref{w-system})
will be presented
in a subsequent article.
The analysis of the accuracy
of the $(m+1)-$st derivative condition
presented in Section~\ref{s-stage1}
carries through,
essentially
(modulo a number of technicalities),
in the more general case as well.
The analysis of the stability
of the functional iteration,
on the other hand,
is far more involved.
The reason for that
is that,
although
the hyperplane $u=u_0$
and
the space tangent
to the fast fibration
over the slow manifold
coincide to leading order
for explicit fast--slow systems (\ref{z-system}),
this is not the case
for the more general systems (\ref{w-system}).
The absence of this feature
makes the stability question
for the functional iteration
far more difficult to answer
in the general case.

In addition,
we are in the process
of generalizing the results
of this article
to other maps
that may be used
in the context
of the functional iteration scheme
developed in \cite{GKKZ-2005}.
In particular,
it is of interest
to use maps
which are implicitly defined
(as opposed to
the explicitly defined ones
presented in \cite{GKKZ-2005}
and in this article).
Preliminary analytical results
for $m=0$ and $m=1$
indicate that
one may construct
functional iteration schemes
based on {\em implicit} maps
which not only retain the accuracy
of the functional iteration scheme
presented in this article
but which are also unconditionally stable.
Moreover, we think that
this analysis may be extended
to higher values of $m$,
and we note that
it is also possible
to carry out the functional iteration
with implicitly defined maps
even when one only has a legacy code
as a timestepper.

\appendix

\section{The one-higher-order proposition
\label{s-1StepAhead}}
\setcounter{equation}{0}
In this appendix,
we state and prove a technical proposition
--
called the one-higher-order proposition
--
about the asymptotic accuracy
of approximations of $\Me$
given an approximation
of the normal space
to $\Me$.
This result is instrumental
in the proof
of the technical lemmas
contained
in the next appendix.

We begin by recalling
the useful formulation,
Eq.~(\ref{inv-eq-alt}),
of the invariance equation
that defines the function $h(x)$,
whose graph is the invariant,
slow manifold $\Me$.
This formulation
revealed that
the matrix
$(-Dh(x),I_{N_f})$
forms a basis
for $\N_z \Me$,
the space normal
to the slow manifold
at the point
$z=(x,h(x)) \in {\mathcal L}$.

The function
$h(x)$
admits an asymptotic expansion
in $\eps$,
\begin{equation}
  h(\cdot)
=
  \sum_{i=0} \eps^i \hi (\cdot) ,
\Label{he-exp}
\end{equation}
where
the coefficients $\hi$,
$i = 0, 1, \ldots\,$,
are determined
by expanding asymptotically
the left member of Eq.~(\ref{inv-eq})
and
setting the coefficient
of $\eps^i$
equal to zero
to obtain
\[
    g_i
    - \sum_{\ell=0}^{i-1} (D\hell) f_{i-1-\ell}
=
    0 ,
    \quad i = 0, 1, \ldots\,,
\]
where the sum is understood to be empty
for $i=0$.
The first few equations are
\begin{eqnarray}
  &&g_0
=
  0 ,
\Label{h(0)} \\
  &&(D_y g)_0 \hone
+
  (D_\eps g)_0
-
  (D\ho) f_0
=
    0 .
\Label{h(1)}
\end{eqnarray}
Here,
Eq.~(\ref{h(0)})
is satisfied identically,
Eq.~(\ref{h(1)})
yields the coefficient $\hone$,
and so on.

The one-higher-order proposition,
which we now state and prove,
establishes
a connection
between
the order in $\eps$
to which
a set $N$
of row vectors
approximates $\N_z \Me$
and
the order to which
the solution $\eta(x)$
to the condition
$N\, G = 0$
approximates $h$.
\begin{PROPOSITION}
\label{l-1stepahead}
Let $N(x,\eps)$ be
an $\Nf \times \Nfs$ matrix
with the property that
its rows span $\N_z\Me$
up to and including
terms of ${\mathcal O}(\eps^m)$,
for some $m = 0, 1, \ldots\,$.
That is,
$N(\cdot,\eps)$
is of the form
\be
  N(\cdot,\eps)
=
  C
  \left(
  - \sum_{i=0}^m \eps^i D\hi(\cdot)
  - \sum_{i \ge m+1} \eps^i R_i(\cdot)
\, , \,
  I_\Nf
  \right) ,
\Label{N=C(-Dh,I)}
\ee
where
$C$ is a non-singular
$\Nf \times \Nf$
matrix
and
$R_i \ne  D\hi$,
for $i = m+1, m+2, \ldots\,$,
in general.
Then,
the condition
\begin{eqnarray}
  N(x, \eps)\, G(x, y, \eps) = 0
\Label{NG=0}
\end{eqnarray}
can be solved for $y$
to yield a function
$y = \eta(x)$,
the asymptotic expansion of which
agrees with that of $h(x)$
up to and including
terms of ${\mathcal O}(\eps^{m+1})$,
\begin{eqnarray}
  \eta(x)
=
  \sum_{i=0} \eps^i \eta_i(x)
=
  \sum_{i=0}^{m+1} \eps^i \hi(x)
+
  {\mathcal O}(\eps^{m+2}) .
\Label{psi=h}
\end{eqnarray}
\end{PROPOSITION}
This proposition is called
the one-higher-order proposition,
because
it states
that the order
to which
$\eta(x)$ approximates
the full slow manifold
is
of one higher
than that
to which
$N$ approximates
the normal space.
\paragraph{Proof of Proposition~\ref{l-1stepahead}.}
We recall
that $h(\cdot)=\Sigma_{i=0} \eps^i \hi (\cdot)$,
by Eq.~(\ref{he-exp}), and
that $\hi$ is determined
from the ${\mathcal O}(\eps^i)$ terms
of the invariance
equation~(\ref{inv-eq-alt}).
Similarly,
$\eta_i$ is determined
from the
${\mathcal O}(\eps^i)$ terms
of Eq.~(\ref{NG=0}).
Thus,
to establish Eq.~(\ref{psi=h}),
it suffices to compare
the terms
of these two equations
from ${\mathcal O}(1)$
up through and including
${\mathcal O}(\eps^{m+1})$
and to show
that they are equal.

First,
for each
$i=0,1,\ldots,m,$
the invariance equation~(\ref{inv-eq-alt})
at ${\mathcal O}(\eps^i)$ is
\be
  \left( - D\ho , I_\Nf \right)
  G_i
+
  \sum_{\ell=1}^i
  \left( - D\hell , 0 \right)
  G_{i-\ell}
=
  0 .
\Label{inv-eq-alt-i}
\ee
Second,
to derive
the ${\mathcal O}(\eps^i)$ terms
for the condition
$NG=0$,
we substitute
the hypothesis~(\ref{N=C(-Dh,I)})
in Eq.~(\ref{NG=0})
and
left-multiply
by $C^{-1}$
to obtain
\be
  C^{-1}\, N\, G
=
  \left( - \sum_{i=0}^m \eps^i D\hi + {\mathcal O}(\eps^{m+1}),
           I_\Nf
  \right)\, G
=
  0 .
\Label{CinvNG=0}
\ee
Thus,
for each
$i=0,1,\ldots,m$,
this condition
at ${\mathcal O}(\eps^i)$ is
\[
    \left( - D\ho , I_\Nf \right)
    G_i
    + \sum_{\ell=1}^i
    \left( - D\hell , 0 \right)
    G_{i-\ell}
=
    0 .
\]
Plainly,
this equation
is identical
to Eq.~(\ref{inv-eq-alt-i}).
Thus,
$\eta_i = \hi$,
for $i = 0, 1, \ldots, m$.

Finally,
we look at
the ${\mathcal O}(\eps^{m+1})$ terms
of the two equations.
Eq.~(\ref{inv-eq-alt-i})
with $i=m+1$ is
\be
  \left( - D\ho , I_\Nf \right)
  G_{m+1}
  + \sum_{\ell=1}^m
  \left( - D\hell , 0 \right)
  G_{m+1-\ell}
  + \left( - D\hmone , 0 \right)
  G_0
=
    0 .
\Label{inv-eq-alt-m+1}
\ee
Also,
Eq.~(\ref{CinvNG=0})
at ${\mathcal O}(\eps^{m+1})$ is
\be
  \left( - D\ho , I_\Nf \right)
  G_{m+1}
  + \sum_{\ell=1}^m
  \left( - D\hell , 0 \right)
  G_{m+1-\ell}
  + \left( R_{m+1} , 0 \right)
  G_0
=
  0 .
\Label{CinvNG-m+1=0}
\ee
We note that $R_{m+1} \ne - D\hmone$,
in general.
However,
$G_0=0$,
since the terms
appearing in
Eqs.~(\ref{inv-eq-alt-m+1})--(\ref{CinvNG-m+1=0})
are evaluated at
$(x, \eta_0, 0) \equiv (x, \ho, 0)$.
Thus,
Eqs.~(\ref{inv-eq-alt-m+1})
and (\ref{CinvNG-m+1=0})
also agree,
and hence
$\eta_{m+1} = \hmone$.
This completes the proof
of the proposition.
\qed

\section{Proofs of Lemmata \ref{l-KGKvsCSP} and \ref{l-cond}
\label{s-derivcond}}
\setcounter{equation}{0}
In this appendix,
we prove lemmata~\ref{l-KGKvsCSP}
and
\ref{l-cond}
characterizing
the asymptotic accuracy
of the approximation to $\Me$
obtained from
the $(m+1)-$st derivative condition
(\ref{cond-alt}).

\paragraph{Proof of Lemma~\ref{l-KGKvsCSP}.}
We write
$z_m$ for $(x, \hm(x))$
and
$z$ for $(x,h(x))$.
The strategy
is as follows:
We will show
that the rows of
$(D_z L_m)(z_m , \eps)$
span $\N_z \Me$
up to and including
terms of ${\mathcal O}(\eps^m)$.
Then,
we will apply
Proposition~\ref{l-1stepahead}
to establish
Eq.~(\ref{lemma51-exp}).

The manifold $\Lm$
is the graph
of the function
$\hm$,
and thus
it coincides exactly
with the zero level set
of the function $- \hm(x) + y$.
As a result,
the rows of
the $\Nf \times \Nfs$
gradient matrix
$(- D\hm(x) , I_\Nf )$
form a basis
for $\N_{z_m} \Lm$.
Second,
the function $\hm(\cdot)$
is defined
through the $(m+1)-$st
derivative condition
$L_m(\cdot, \hm(\cdot), \eps) = 0$.
Therefore,
$\Lm$ also coincides with
(a connected component of)
the zero level set
of the function
$L_m(z, \eps)$.
Thus,
the rows
of the $\Nf \times \Nfs$
gradient matrix
$(D_z L_m)(x, \hm(x), \eps)$
also form a basis
for $\N_{z_m} \Lm$.
It follows
from the existence
of these two bases
that there exists
a non-singular $\Nf \times \Nf$
matrix $C$
such that
\begin{equation}
  \left(D_z L_m\right)(\cdot, \hm(\cdot), \eps)
=
  C \left(- D\hm(\cdot) , I_\Nf \right) .
\Label{DL-vs-Dhm}
\end{equation}

Next,
the induction hypothesis implies
that the asymptotic expansions
of $\hm$ and $h$
agree up to and including
terms of ${\mathcal O}(\eps^m)$,
\begin{equation}
  \hm(\cdot)
=
  \sum_{i=0}^m \eps^i \hi(\cdot)
+
  {\mathcal O}(\eps^{m+1}) .
\Label{hm-vs-h}
\end{equation}
Since
the vector field
is assumed to be
sufficiently smooth,
we may differentiate
both sides of this equation
with respect to $x$
to obtain
\begin{equation}
  D\hm(\cdot)
=
  \sum_{i=0}^m \eps^i D\hi(\cdot)
+
  {\mathcal O}(\eps^{m+1}) .
\Label{Dhm-vs-Dh}
\end{equation}
Combining Eqs.~(\ref{DL-vs-Dhm})
and
(\ref{Dhm-vs-Dh}),
then,
we find
\[
  \left(D_z L_m\right)(\cdot, \hm(\cdot), \eps)
=
  C
\left(
  - \sum_{i=0}^m \eps^i D\hi(\cdot) + {\mathcal O}(\eps^{m+1}) ,
  I_\Nf
\right) .
\]
This equation shows
that the rows
of $(D_z L_m)(x, \hm(x), \eps)$
span $\N_z \Me$
up to and including
terms of
${\mathcal O}(\eps^m)$.
Hence,
application
of the one-higher-order proposition,
Proposition~A.1,
completes the proof of this lemma.
\qed
Before we proceed
with the proof
of Lemma~\ref{l-cond},
we prove
the following result
which
will be needed therein.
\begin{LEMMA}
\label{l-L(m)}
For $m=0, 1, \ldots$,
for $H = {\mathcal O}(\eps)$,
and
for a general point
$z = (x,y)$,
the function $L_m$
is written as
\[
  L_m(z)
=
  ( -\eps^{-1} H )^{m+1}
\left[ \left(D_y g\right)_0(z) \right]^m
  g_0(z)
+
  {\mathcal O}\left( \eps, \|g_0(z)\|^2 \right) ,
\]
where the notation ``$(\cdot)_0(z)$''
stands for $(\cdot)(z , 0)$.
The Jacobian $D_yL_m$
is written as
\be
  \left(D_y L_m\right)(z)
=
  \left( -\eps^{-1} H \left(D_y g\right)_0 \right)^{m+1}
+
  {\mathcal O}\left(\eps, \|g_0(z)\| \right) .
\Label{DyLm-O(1)}
\ee
\end{LEMMA}
\begin{PROOF}
For this proof,
we write
$(\cdot)_0$
instead of
$(\cdot)_0(z)$
for the sake of brevity.
The proof
is by induction
on $m$.
For $m=0$,
we recall
Eq.~(\ref{L(0)-raw}),
\begin{eqnarray*}
  L_0
&=&
  -\eps^{-1} H g ,
\end{eqnarray*}
and hence,
expanding $g$
in powers of $\eps$,
we find
\begin{eqnarray*}
  L_0
=
  -\eps^{-1} H g_0
+
  {\mathcal O}(\eps) .
\end{eqnarray*}
This is
the desired formula
for $L_0$.
Differentiating
both members
of this formula
with respect to $y$,
we obtain
\[
  D_y L_0
=
  -\eps^{-1} H \left(D_y g\right)_0
+
  {\mathcal O}(\eps) .
\]
This is
the desired formula
for $D_y L_0$.

Next,
we carry out
the induction step
for general $m$,
namely
we assume that
\begin{eqnarray}
  L_m
&=&
  \left( -\eps^{-1} H \right)^{m+1}
  \left(D_y g \right)_0^m
  g_0
+
  {\mathcal O}\left( \eps , \|g_0\|^2 \right) ,
\Label{ind-hypo-Lm}
\\
  D_y L_m
&=&
  \left( -\eps^{-1} H \left(D_y g\right)_0 \right)^{m+1}
+
  {\mathcal O}\left(\eps, \|g_0(z)\| \right)
\Label{ind-hypo-DyLm}
\end{eqnarray}
and
show that
\begin{eqnarray}
  L_{m+1}
&=&
  \left( -\eps^{-1} H \right)^{m+2}
  \left(D_y g \right)_0^{m+1}
  g_0
+
  {\mathcal O}\left( \eps , \|g_0\|^2 \right) .
\Label{ind-step-Lm}
\\
  D_y L_{m+1}
&=&
  \left( -\eps^{-1} H \left(D_y g\right)_0 \right)^{m+2}
+
  {\mathcal O}\left(\eps, \|g_0(z)\| \right) .
\Label{ind-step-DyLm}
\end{eqnarray}
By Eq.~(\ref{L(m+1)=DL(m)G}),
\[
  L_{m+1}
=
  -\eps^{-1} H
   (D_z L_m) G
=
  -\eps^{-1} H
\left[
  \eps (D_x L_m) f
+
  (D_y L_m) g
\right] ,
\]
Then,
we substitute
the induction
hypothesis~(\ref{ind-hypo-Lm})
into this expression.
Application
of the differential operator
$(-H/\eps)[\eps (D_x \cdot) f + (D_y \cdot)g]$
on the ${\mathcal O}( \eps , \|g_0\|^2 )$
remainder
does not alter
its asymptotic magnitude.
Moreover,
the term
$\eps (D_x L_m) f$
is ${\mathcal O}(\eps)$
and, hence,
can be absorbed
also in the remainder.
Therefore,
we are left with
the term $(-H / \eps)(D_yL_m) g$.
Substituting
$D_yL_m$
into this expression
from the induction hypothesis~(\ref{ind-hypo-DyLm}),
we arrive
at the desired formula~(\ref{ind-step-Lm}).

Finally,
we prove
the leading order formula~(\ref{ind-step-DyLm}).
First,
we differentiate
both members of
the leading order formula~(\ref{ind-step-Lm})
with respect to $y$
and
use the product rule derivative
to evaluate
the right member.
The second term
from the product rule
is precisely the leading order term
in Eq.~(\ref{DyLm-O(1)}).
The other term
from the product rule,
\[
  m
  \left( -\eps^{-1} H \right)^{m+2}
  \left(D^2_y g\right)_0
  \left( (D_y g)_0^{m-1}, g_0\right) ,
\]
may be absorbed
in the remainder
since it is
linear in $g_0$.
Thus,
we have obtained
the desired
formula~(\ref{ind-step-DyLm})
and completed the proof
of the lemma.
\end{PROOF}
\paragraph{Proof of Lemma~\ref{l-cond}.}
We first use \cite[Theorem 3]{C-1981}
to establish
that condition~(\ref{cond-alt})
has a solution $y = \hmm(x)$
which is
${\mathcal O}(\eps^{m+1})-$close
to $\tilde{h}_{m+1}$.
According to that theorem,
it suffices to show that
\[
  \left((D_z L_m)(x, \tilde{h}_{m+1}(x), \eps)\right)
  G(x, \tilde{h}_{m+1}(x), \eps)
=
  {\mathcal O}(\eps^{m+2}) .
\]
By the definition
of $\tilde{h}_{m+1}$,
\[
  \left((D_z L_m)(\cdot, \hm(\cdot), \eps)\right)
  G(\cdot, \tilde{h}_{m+1}(\cdot), \eps)
=
  0 .
\]
Thus,
we may write
\begin{eqnarray}
\lefteqn{
  \left((D_z L_m)(\cdot, \tilde{h}_{m+1}(\cdot), \eps)\right)
  G(\cdot, \tilde{h}_{m+1}(\cdot), \eps)
}
\nonumber\\
&=&
  \left[
      (D_z L_m)(\cdot, \tilde{h}_{m+1}(\cdot), \eps)
    - (D_z L_m)(\cdot, \hm(\cdot), \eps)
  \right]
  G(\cdot, \tilde{h}_{m+1}(\cdot), \eps) .\hspace{0.8cm}
\Label{DLg-aux}
\end{eqnarray}

Next,
we have the following
estimates
of the asymptotic magnitudes
of the two terms
in the right member
of Eq.~(\ref{DLg-aux}):
\[
  \tilde{h}_{m+1}
=
  \sum_{i=0}^{m+1} \eps^i \hi + {\mathcal O}(\eps^{m+2})
\]
by Lemma~\ref{l-KGKvsCSP},
and also
\[
  \hm
=
  \sum_{i=0}^m \eps^i \hi + {\mathcal O}(\eps^{m+1})
\]
by the induction hypothesis.
Thus,
\begin{eqnarray*}
    \tilde{h}_{m+1}
-
    \hm
=
    {\mathcal O}(\eps^{m+1}) ,
\end{eqnarray*}
and hence
Taylor's Theorem with remainder
yields
\begin{eqnarray}
  (D_z L_m)(\cdot, \tilde{h}_m+1(\cdot), \eps)
-
  (D_z L_m)(\cdot, \hm(\cdot), \eps)
=
  {\mathcal O}(\eps^{m+1}) ,
\Label{est-first-term}
\end{eqnarray}
since $L_m$
and its derivatives
are ${\mathcal O}(1)$.
This is
the desired estimate
of the first term
in the right member
of Eq.~(\ref{DLg-aux}).

It remains to estimate
the second term,
$G(\cdot, \tilde{h}_{m+1}(\cdot), \eps)$
in the right member
of Eq.~(\ref{DLg-aux}).
We recall that
$G= \left( { \eps f \atop g } \right)$,
where $f$ and $g$
are ${\mathcal O}(1)$
in general.
Hence,
the first component
of $G(\cdot, \tilde{h}_{m+1}(\cdot), \eps)$
is plainly ${\mathcal O}(\eps)$.
The second component
is as well,
since
Lemma~\ref{l-KGKvsCSP}
implies that
$\tilde{h}_{m+1,0} = \ho$
and hence that
$g(\cdot, {\tilde h}_{m+1} (\cdot), \eps)
= {\mathcal O}(\eps)$, also.
Therefore,
\be
  G(\cdot, \tilde{h}_{m+1}(\cdot), \eps)
=
  {\mathcal O}(\eps).
\Label{est-second-term}
\ee

Combining the estimates
(\ref{est-first-term})
and
(\ref{est-second-term}),
we see
that the right member
of Eq.~(\ref{DLg-aux})
is ${\mathcal O}(\eps^{m+2})$,
which is
the desired result.

Finally,
the solution
of the condition
$L_{m+1} = 0$
yields
an $\Ns-$dimensional
manifold $\Lmm$,
as may be shown
using the Implicit Function Theorem
and
\cite[Theorem~1.13]{O-1986}.
It suffices
to show that
\[
  \mathrm{det}\left(D_y L_{m+1} \right)(\cdot, \hmm(\cdot)) \ne 0 .
\]
Lemma~\ref{l-L(m)}
yields
a leading order formula
for $D_yL_{m+1}$,
\[
  \left(D_y L_{m+1}\right)(z)
=
  \left( -\eps^{-1} H \left(D_y g\right)_0 \right)^{m+2}
+
  {\mathcal O}\left(\eps, \|g_0(z)\| \right) .
\]
Here,
$z$ is a general point
and
$(\cdot)_0(z) = (\cdot)(z , 0)$.
Next,
we showed above that
$h_{(m+1,0)} = h_0$.
Recalling, then,
Eq.~(\ref{inveq-0}),
we obtain
\[
  \left(D_y L_{m+1}\right)(x, \hmm(x))
=
  \left[ -\eps^{-1} H \left(D_y g\right)_0 \right]^{m+2}
+
  {\mathcal O}(\eps) ,
\quad \mbox{for all} \
  x \in K,
\]
where
$(D_y g)_0
=
(D_y g)(x , h_0(x) , 0)$.
Thus,
\[
  \mathrm{det}\left(D_y L_{m+1} \right)(x, \hmm(x))
\ne
  0 ,
\quad \mbox{for all} \
  x \in K,
\]
by normal hyperbolicity
and
the proof is complete.
\qed

\end{document}